\documentclass[11pt,a4paper]{article}
\usepackage[active]{srcltx}
\usepackage{amsmath,amssymb,amsthm,amsfonts}
\usepackage{pictex}
\usepackage{comment}
\usepackage{rotating}
\usepackage{color}

\usepackage{graphics}
\numberwithin{equation}{section}
\def\supp{\textrm{supp}}

\def\~{\widetilde}\def\hat{\widehat}\def\g{{\geqslant}}
\def\R{\operatorname{Re}}
\def\Ima{\operatorname{Im}}

\newcommand{\Pp}{\mathcal P}\newcommand{\Ee}{\mathcal E}
\newcommand{\RR}{\mathbb R}\newcommand{\Rr}{\mathcal R}
\newcommand{\CC}{\mathbb C}

\newcommand{\refe}[1]{(\ref{#1})}

\newtheorem{definition}{Definition}[section]
\newtheorem{theorem}{Theorem}[section]
\newtheorem{example}{Example}[section]
\newtheorem{lemma}{Lemma}[section]
\newtheorem{proposition}{Proposition}[section]

\newtheorem{remark}[theorem]{Remark}

\newcommand{\abs}[1]{\left\vert#1\right\vert}

\begin{document}

\title{Asymptotics for Multiple Meixner Polynomials}

\author{A. Aptekarev \thanks{The research of A. Aptekarev was supported by the grants
RFBR 11-01-12045 OFIM,  11-01-00245 and the Chair Excellence Program of
Universidad Carlos III Madrid, Spain and Bank Santander. }, and J.
Arves\'u\thanks{The research of J. Arves\'u was partially supported by the
research grant MTM2009-12740-C03-01 of the Ministerio de Educaci\'on y Ciencia
of Spain and grants CCG07-UC3M/ESP-3339 and CC\-G08-UC3M/ESP-4516 from
Comunidad
Aut\'onoma de Madrid.} \medskip \\
Keldysh Institute for Applied Mathematics,
Russian Academy of Sciences\\
Miusskaya pl. 4, 125047 Moscow,
RUSSIA\\
Department of Mathematics, Universidad Carlos III de Madrid,\\
Avda. de la Universidad, 30, 28911, Legan\'es, Madrid, SPAIN\\
}

\maketitle

\begin{abstract}
We study the asymptotic behavior of Multiple Meixner polynomials of first and
second kind, respectively \cite{arvesu_vanAssche}. We use an algebraic function
formulation for the solution of the equilibrium problem with constrain to
describe  their zero distribution. Then analyzing  the limiting behavior  of
the coefficients of the recurrence relations for Multiple Meixner polynomials
we obtain the main term of their asymptotics.
\end{abstract}



\section{Introduction\label{introduction}}

We consider the asymptotics of polynomial sequences $(P_n)$ defined by
orthogonality relations with respect to  a discrete measure $\mu$ (with finite
moments)
\begin{equation}\label{discrmes}
\mu=\sum_{k=0}^{N}\rho(x)\delta_{x_k},\quad \rho(x_k)>0,\ x_k\in\mathbb{R}
,\quad N\in\mathbb{N}\cup \{+\infty\},
\end{equation}
which is a linear combination of Dirac measures at the points $x_0,\ldots,x_N
$. By $\mathbb{N}$ we denotes the set of all nonnegative integers.

The asymptotic theory of orthogonal polynomials with respect to discrete measures is not so widely developed as for absolutely continues measures of orthogonality. It can be explained because the distribution of zeros of the polynomials orthogonal with respect to a discrete measure has an extra constrain. Indeed (it follows from the interlacing property of orthogonal polynomials), between two neighboring mass points of the measure of orthogonality it can not be situated more than one zero of the orthogonal polynomials. This constrain plays an important role in the logarithmic potential description of the limiting measure of zeros distribution and of the main term of asymptotics.

Logarithmic potential
$$\Pp^{\nu}(z)\,=\,-\int\ln|z-t|\,\,d\nu(t),$$
for probability measure $\nu_{P_n}(t)$ (\textit{zero counting measure}), equally
distributed in the zeros of polynomial $P_n$, is
$$
\Pp^{\nu_{P_n}}\,=\,-\frac1n\ln|P_n|\,.
$$
Thus, the potential of a limiting (when $n\!\to\!\infty$) zero distribution
measure $\lambda$ defines an \textit{exponent of the main term} of the
polynomial sequence $(P_n)$ asymptotics. Such type of asymptotics are
also called \textit{weak asymptotics}. For the description of this limiting measure it is useful to formulate a minimization problem for the energy of the logarithmic potential. It is clear, that constrain of the zero counting measure (which we pointed out above) leads to a problem of energy minimization in the  class of probabilistic measures constrained by the weak limit of the mass points counting measure. For discrete orthogonal polynomials this approach was suggested by Rakhmanov in \cite{RaDiscr1}. This approach was later extended in \cite{SD}
(for more details on the topic we refer to monographs \cite{Nikishin} and \cite{SafTot} for extremal and equilibrium problem of logarithmic potential with
application to the weak asymptotics of polynomials sequences as well as to
monograph \cite{BKMM} for discrete orthogonality and constrained equilibrium).

The content of this paper starts with the study of weak asymptotics for monic polynomial sequences defined by discrete measure \eqref{discrmes}. In particular, we focus our attention in a polynomial sequence orthogonal with respect to negative binomial distribution (Pascal distribution) on $\mathbb{N}$, i.e
\begin{equation}\label{meixmes}
\rho_k=\frac{(\beta)_k}{k!}c^k, \quad (\beta>0\,,\,\,0<c<1) \quad
\mbox{on}\quad x_k=k,\quad k\in \mathbb{N}.\end{equation} This orthogonal
polynomial sequence (denoted by $M_n(x;\beta,c)$) is called \textit{classical
Meixner polynomials} \cite{Niki}. They satisfy the orthogonality conditions
\begin{equation}\label{meixpol}\sum \limits_{k=0}^{+\infty}M_n(k;\beta,c)(-k)_j\frac{(\beta)_k}{k!}c^k
=0, \qquad j=0,1,\ldots,n-1,\end{equation} where $(-k)_{j}=-k(1-k)\cdots
(j-k-1)$, $j\in\mathbb{N}$, and $(k)_{0}=1$, is the Pochhammer symbol. However,
our main goal is the asymptotic analysis of \textit{multiple Meixner
polynomials}, which constitute a generalization of the aforementioned Meixner
polynomials. For these polynomials the orthogonality conditions are considered
with respect to a collection of Meixner weights \eqref{meixmes} with different
parameters $\beta$ or $c$. The study of these polynomials (among others
classical discrete multiple orthogonal polynomials) was initiated in
\cite{arvesu_vanAssche}.

The structure of the paper is as follows. The next subsection~\ref{Ss1.1} is
devoted to the notion of multiple orthogonal polynomials (definition and some
properties). In subsection~\ref{recursiemeervoudige} we give general
information about recurrence relations for multiple orthogonal polynomials. The
introduction finishes with two subsections devoted to some formal properties of
the multiple Meixner orthogonal polynomials obtained in
\cite{arvesu_vanAssche}. Recall that there are two kind of multiple Meixner
orthogonal polynomials. The first kind corresponds (see subsection~\ref{mmvb1})
to simultaneous orthogonality conditions with respect to discrete weights
formed by \eqref{meixmes} with various $c_i, \,\, i=1, \ldots, r$. The second
kind polynomials (see subsection~\ref{mmvb2}) appear when the collection of
discrete weights is formed by \eqref{meixmes} with various $\beta_i, \,\,
i=1,\ldots, r$.  The explicit expressions for the coefficients of the
recurrence relations of the multiple Meixner orthogonal polynomials are the
main outcome from these subsections~\ref{mmvb1} and~\ref{mmvb2} -used in the
sequel.

In section~\ref{equil-alg} we state the logarithmic potential equilibrium problems for the description of weak asymptotics of Meixner and multiple Meixner polynomials and the solutions of these problems by means
of certain algebraic functions are found. The procedure for solving such type of equilibrium problems is inspired in \cite{AKLT} (see also earlier references therein).

Lastly, in section~\ref{nroot} we, starting from the coefficient
of the recurrence relations, obtain the main term of asymptotics of multiple
Meixner polynomials and then we check the connection of this term with the equilibrium problem from the previous section. Theorems~\ref{T3} and~\ref{T3.2} are the main results of the paper. 

Furthermore, the spectral curves for the multiple Meixner orthogonal polynomials (of the first and second kind) and their connection with recurrence relations and vector potential equilibrium with external field and constrain constitute the main outcome of the present paper.

Finally, we should also mention paper \cite{Sorokin} by Sorokin in which the weak asymptotics of multiple Meixner orthogonal polynomials of the second type were studied, but using another approach where the integral representation for the polynomials and saddle point method were applied.

\subsection{Multiple orthogonal polynomials (type II)}\label{Ss1.1}
Below we summarize some results needed for the sequel.

Assume one has $r$ positive measures $\mu_1,\dots,\mu_r$ on $\mathbb{R}$,
where for each $i=1,2,\dots,r$
\[\supp(\mu_i)=\left\{x\in\mathbb{R}\left|\
\forall \epsilon>0,\
\mu_i\Bigl((x-\epsilon,x+\epsilon)\Bigr)>0\right.\right\},\]
denotes the support of $\mu_i$. By $\Omega_i$ we denote the smallest interval that contains $\supp(\mu_i)$.

We introduce a multi-index
$\vec{n}=(n_1,n_2,\ldots,n_r)\in\mathbb{N}^r$ and its length $|\vec{n}|=n_1+n_2+\ldots+n_r$. Multiple orthogonal polynomials (of type II)
can be defined as follows \cite[Chapter 4.3]{Nikishin}, \cite{aptekarev} and
\cite{artikel}.

\begin{definition}
A type II multiple orthogonal polynomial $P_{\vec{n}}$, corresponding to the
multi-index $\vec{n}\in\mathbb{N}^r$, is a polynomial of degree $\le |\vec{n}|$ which satisfies the orthogonality conditions
\begin{equation*}
\int_{\Omega_i}P_{\vec{n}}(x)x^k\ d\mu_i(x) =  0,\quad k=0,1,\ldots,n_i-1,\quad
i=1,2,\dots,r.
\label{stelseltype2}
\end{equation*}
\end{definition}

Here one has a linear system of $|\vec{n}|$ homogeneous relations for the $|\vec{n}|+1$ unknown coefficients of $P_{\vec{n}}$.  Notice that the solution $P_{\vec{n}}$ which is unique up to a multiplicative factor and also that this polynomial solution has exactly degree $|\vec{n}|$ (then the monic multiple orthogonal polynomial exists and will be unique).  When this happens $\vec{n}$ is said to be a \textit{normal index} for type II.

Suppose that the above $r$ positive measures on $\mathbb{R}$ are discrete, i.e.
\[\mu_i=\sum_{k=0}^{N_i}\rho_{i,k}\delta_{x_{i,k}},\qquad \rho_{i,k}>0,\ x_{i,k}\in\mathbb{R},\ \
N_i\in\mathbb{N}\cup \{+\infty\},\ \  i=1,\ldots,r,\]
with all the $x_{i,k}$,
$k=0,\ldots,N_i$, different and this for each $i$.  In this case we have that $\supp(\mu_i)$ is the closure of  $\{x_{i,k}\}_{k=0}^{N_i}$ and that $\Omega_i$ is the smallest closed interval on $\mathbb{R}$ which contains $\{x_{i,k}\}_{k=0}^{N_i}$. The corresponding polynomials are then discrete multiple orthogonal polynomials \cite{arvesu_vanAssche}. The discrete measures in this paper will be supported on $\mathbb{N}$ or a subset, which is
achieved by taking $x_{i,k} = k$ for $i=1,\ldots,r$. The orthogonality
conditions are then more conveniently expressed in terms of the
polynomials $(-x)_j$.

\begin{definition}
A discrete multiple orthogonal polynomial of type II on the linear lattice, corresponding to the multi-index
$\vec{n}\in\mathbb{N}^r$, is a polynomial
 $P_{\vec{n}}$ of degree
$\le |\vec{n}|$ that satisfies the orthogonality conditions
\begin{equation*}
\sum \limits_{k=0}^{N_i}P_{\vec{n}}(k)(-k)_j\,\rho_{i,k}=0,\quad j=0,1,\ldots,n_i-1,\quad i=1,2,\dots,r.
\label{steltype2}
\end{equation*}
\end{definition}

In what follows we will focus on a systems of measures
\[\mu_i=\sum_{k=0}^{N} \rho_{i,k} \delta_{x_{k}},\qquad \rho_{i,k}>0,\
x_k\in\mathbb{R},\ \ N\in
\mathbb{N} \cup \{+\infty\},\ \ \ i=1,\ldots,r,
\]
where $supp(\mu_i)$ is the closure of $\{x_k\}_{k=0}^N$ and
$\Omega_i=\Omega$ for each $i=1,\ldots,r$, for which every multi-index is normal, namely an {\it AT system} (see Nikishin and Sorokin \cite[p.~140]{Nikishin} as well as Arves\'u et al. \cite[p.~26]{arvesu_vanAssche}). An interesting property related to the location of the zeros of the discrete multiple orthogonal polynomials is given by the following theorem.

\begin{theorem} {\bf(see \cite[p.~26]{arvesu_vanAssche})}
Suppose we have an AT system of $r$ positive discrete measures. Then every discrete multiple orthogonal polynomial $P_{\vec{n}}$ of type II, corresponding to the multi-index
$\vec{n}$ with $|\vec{n}|< N+1$, has exactly $|\vec{n}|$ different zeros on $\Omega$.
\end{theorem}

Observe that in an AT system every multi-index $\vec{n}$, with $|\vec{n}|< N+1$, is normal. We refer to \cite[p.~26]{arvesu_vanAssche} for the following two examples of Chebyshev systems:
\begin{example}\cite[Example 3 on p.~138]{Nikishin}
\label{example}
The functions
\begin{displaymath}
\begin{array}{c}
v(x)c_1^x,v(x)xc_1^x,\ldots,v(x)x^{n_1-1}c_1^x,\\
\vdots \\
v(x)c_r^x,v(x)xc_r^x,\ldots,v(x)x^{n_r-1}c_r^x,
\end{array}
\end{displaymath}
with all the $c_i>0$, $i=1,\ldots,r$, different and $v$ a function which has no zeros on $\mathbb{R}^+$, form a Chebyshev system on $\mathbb{R}^+$ for every $\vec{n}=(n_1,\ldots,n_r)\in \mathbb{N}^r$.
\end{example}

\begin{example}  \label{example2}
The functions
\begin{equation*}\label{defuncties}
 \begin{array}{c}
v(x) \Gamma(x+\beta_1),v(x)x \Gamma(x+\beta_1),\ldots,v(x)x^{n_1-1}\Gamma(x+\beta_1),\\
\vdots \\
v(x)\Gamma(x+\beta_r),v(x)x\Gamma(x+\beta_r),
\ldots,v(x)x^{n_r-1}\Gamma(x+\beta_r),
\end{array}
\end{equation*}
with $\beta_i >0$ and
$\beta_i - \beta_j \notin \mathbb{Z}$ whenever $i \neq j$ and
$v$ a function with no zeros on $\mathbb{R}^+$, form a Chebyshev
system on $\mathbb{R}^+$ for every $\vec{n} = (n_1,\ldots,n_r) \in
\mathbb{N}^r$. If $\beta_i - \beta_j \notin \{0,1,\ldots,N-1\}$ then this still
gives a Chebyshev system for every $\vec{n}=(n_1,\ldots,n_r)$ for which $n_i < N+1$, $i=1,2,\ldots,r$.
\end{example}

\subsection{Recurrence relation}
\label{recursiemeervoudige}

Suppose that all the multi-indices are normal for the $r$ measures $\mu_1,\ldots, \mu_r$.  Then there exists an interesting recurrence relation of order $r+1$ for the monic multiple orthogonal polynomials of type II with nearly diagonal multi-indices \cite{artikel}.  Here the nearly diagonal multi-index, corresponding to $n$, is given by
\[\vec{s}(n)=(\underbrace{k+1,k+1,\ldots,k+1}_{s\  \textrm{times}}
,\underbrace{k,k,\ldots,k}_{r-s\ \textrm{times}}),
\]
with $n=kr+s$, $0\le s <r$.
If we write
$P_n(x)=P_{\vec{s}(n)}(x)$,
then the following recurrence relation holds:
\begin{equation}
\label{recursie}
xP_n(x) = P_{n+1}(x)+\sum_{j=0}^r a_{n,j}P_{n-j}(x),
\end{equation}
with initial conditions
$P_0=1$ and $P_j=0$, $j=-1,-2,\ldots,-r$.

A more general form of this recurrence relation is given by
\begin{eqnarray*}  \label{recursiegen}
xP_{\vec{n}}(x) & = & P_{\vec{n}+\vec{e}_1}(x)+b_{\vec{n},0}P_{\vec{n}}(x)+\sum_{j=1}^r b_{\vec{n},j}P_{\vec{n}-\vec{v}_j}(x),
\end{eqnarray*}
where $\vec{e}_i$ is the $i$th standard unit vector in $\mathbb{R}^r$ and $\vec{v}_j=\sum\limits_{k=0}^{j-1}\vec{e}_{r-k}$.  In the case $r=2$ the recurrence relation
for the polynomials with nearly diagonal multi-indices (\ref{recursie}) gives the relations
$(n+1,n)\rightarrow (n,n) \rightarrow (n,n-1)\rightarrow (n-1,n-1)$ and $(n+1,n+1)\rightarrow
(n+1,n) \rightarrow (n,n) \rightarrow (n,n-1)$.  The first relation follows from the general
case $(n_1+1,n_2)\rightarrow (n_1,n_2) \rightarrow (n_1,n_2-1) \rightarrow (n_1-1,n_2-1)$ by setting $n_1=n_2=n$.  To obtain the second one we set $n_1=n$ and $n_2=n+1$ and interchange the measures $\mu_1$ and $\mu_2$.

\subsection{Multiple Meixner polynomials (first kind)}
\label{mmvb1}
Let $\mu_1,\ldots,\mu_r$ be defined as
\begin{equation}
\label{mesMeixI} \mu_i=\sum_{k=0}^{+\infty} \frac{(\beta)_kc_i^k}{k!}
\delta_{k},\qquad 0<c_i<1,\quad i=1,\ldots,r,\end{equation}
 with all the parameters $c_i$ different.  The support of these measures is $\mathbb{N}$ and
we have $\Delta_1=\ldots=\Delta_r=\mathbb{R}^+$. If
$\beta\notin\mathbb{N}\setminus\{0\}$, we define the functions $w_i^{\beta}$,
$i=1,\ldots,r$, in $C^{\infty}\left(\mathbb{R}\setminus \{-\beta,-\beta
-1,-\beta-2,\ldots\}\right)$ as
\begin{equation}w_i^{\beta}(x)=\left\{
 \begin{array}{ll}
 \displaystyle{\frac{\Gamma(\beta+x)}{\Gamma(\beta)}
\frac{c_i^x}{\Gamma(x+1)}} & \mbox{if $x \in \mathbb{R}\setminus \left( \{-1,-2,-3,\ldots\}\cup\{-\beta,-\beta -1,-\beta-2,\ldots\}\right)$},\\
  0 & \mbox{if $x \in \{-1,-2,-3,\ldots\}$},
 \end{array} \right.\label{weights-kind-1}
\end{equation}
with simple poles at $-\beta,-\beta -1,-\beta -2,\ldots$.  If $\beta \in \mathbb{N} \setminus \{ 0 \}$ we define
\[w_i^{\beta}(x)=\frac{(x+1)_{\beta-1}}{(\beta-1)!}c_i^x,
\]
which are functions in $C^{\infty}$. These functions can be written as
$w_i^{\beta}(k)=(\beta)_kc_i^k/k!$, $k\in\mathbb{N}$, $i=1,\ldots,r$. By
Example \ref{example} we know that the measures $\mu_1,\ldots,\mu_r$ form an AT
system which guarantees that every multi-index $\vec{n}=(n_1,\ldots,n_r)$ is
normal for these measures.

The monic multiple Meixner polynomials of first kind, corresponding to the multi-index
$\vec{n}=(n_1,\ldots,n_r)$ and the parameters $\beta$, $\vec{c}=(c_1,\ldots,c_r)$, is the unique monic polynomial
$M_{\vec{n}}^{\beta;\vec{c}}$ of degree $|\vec{n}|$ which satisfies the orthogonality conditions
\begin{equation*}
\sum_{k=0}^{+\infty}M_{\vec{n}}^{\beta;\vec{c}}(k)(-k)_jw_i^{\beta}(k)=0,
\qquad j=0,\ldots,n_i-1,\quad i=1,2,\dots,r.
\end{equation*}
In \cite{arvesu_vanAssche} the following Rodrigues formula was found
\begin{equation*}\label{MIRod}
M_{\vec{n}}^{\beta;\vec{c}}(x)=\prod_{k=1}^{r}
\left(\frac{c_k}{c_k-1}\right)^{n_k}
(\beta)_{|\vec{n}|}\frac{\Gamma(\beta)\Gamma(x+1)}{\Gamma(\beta+x)}\prod_{i=1}^{r}\left(\frac{1}{c_i^x}\nabla^{n_i}c_i^x\right)
 \frac{\Gamma(\beta+|\vec{n}|+x)}{\Gamma(\beta+|\vec{n}|)\Gamma(x+1)} .
\end{equation*}
In particular, for $r=2$, some straightforward calculations yield
\begin{eqnarray*}
M_{n_1,n_2}^{\beta;c_1,c_2}(x) & = & \left(\frac{c_1}{c_1-1}\right)^{n_1}
\left(\frac{c_2}{c_2-1}\right)^{n_2}(\beta)_{n_1+n_2}\  F_1
\left(
\begin{array}{c}
-x,-n_1,-n_2;\beta;1-\frac{1}{c_1},1-\frac{1}{c_2}
\end{array}
\right)\\
& = &
\frac{c_1^{n_1}c_2^{n_2}(\beta)_{n_1+n_2}
}{(c_1-1)^{n_1}(c_2-1)^{n_2}
}\sum_{j=0}^{n_1+n_2}\sum_{k=0}^{j}
\frac{(-n_1)_k(-n_2)_{j-k}}{(\beta)_{j}}
\frac{\left(\frac{c_1-1}{c_1}\right)^k}{k!}
\frac{\left(\frac{c_2-1}{c_2}\right)^{j-k}}{(j-k)!}(-x)_{j},
\end{eqnarray*}
where
\begin{eqnarray*}
F_1(\alpha,\beta,\beta';\gamma;x,y)
& = &
\sum_{m=0}^{+\infty}\sum_{n=0}^{+\infty}
\frac{(\alpha)_{m+n}(\beta)_m(\beta')_n}{(\gamma)_{m+n}m!n!}x^my^n,
\end{eqnarray*}
is the first of Appell's hypergeometric functions of two variables
\cite{hyperfuncties}. From this explicit expression can be found  the
coefficients of the recurrence relation.
\begin{theorem}\label{1.1}{\bf(see \cite{arvesu_vanAssche})}
Multiple Meixner polynomials (first kind) $M_{n_1,n_2}^{\beta; \,c_1,c_2}$
satisfy the following recurrence relations
\[xM^{\beta;\,c_1,c_2}_{n_1,n_2}(x)=M^{\beta;\,c_1,c_2}_{n_1+1,n_2}(x)+
b_{n_1,n_2}M^{\beta;\,c_1,c_2}_{n_1,n_2}(x)+c_{n_1,n_2}M^{\beta;\,c_1,c_2}_{n_1,n_2-1}(x)+
d_{n_1,n_2}M^{\beta;\,c_1,c_2}_{n_1-1,n_2-1}(x),
\]
where coefficients are
\begin{eqnarray}
b_{n_1,n_2} & = &
n_1(2a_1+1)+n_2(a_1+a_2+1)+a_1\beta,  \nonumber
\\  c_{n_1,n_2}  & = & \Bigl(n_1(a_1^2+a_1)+n_2(a_2^2+a_2)\Bigr)(n_1+n_2+\beta-1), \label{MIrec}
 \\
 d_{n_1,n_2} & = &
 (\beta+n_1+n_2-1)(\beta+n_1+n_2-2)(a_1+1)(a_1-a_2)a_1n_1.  \nonumber
\end{eqnarray}
with $a_1=\frac{c_1}{1-c_1}$ and $a_2=\frac{c_2}{1-c_2}$. \end{theorem}

\subsection{Multiple Meixner polynomials (second kind)}\label{mmvb2}

Assume that the measures $\mu_1,\ldots,\mu_2$ are given by
\begin{equation}\label{mes-2-kind}
\mu_i=\sum_{k=0}^{+\infty} \frac{(\beta_i)_kc^k}{k!} \delta_{k},\qquad
\beta_i>0,\quad i=1,\ldots,r,\end{equation} with $0<c<1$ and all the $\beta_i$
different.  The support of these measures is $\mathbb{N}$.  In
\cite{arvesu_vanAssche} was found that every multi-index
$\vec{n}=(n_1,\ldots,n_r)$ is normal for the above system of measures
$\mu_1,\ldots,\mu_r$, whenever \textbf{$\beta_i-\beta_j \notin \mathbb{Z}$} for
all $i\not=j$.

Here
\begin{equation}
w^{\beta_i}(x)=\left\{
 \begin{array}{ll}
 \displaystyle{\frac{\Gamma(\beta_i+x)}{\Gamma(\beta)}
\frac{c^x}{\Gamma(x+1)}} & \mbox{if $x \in \mathbb{R}\setminus \left( \{-1,-2,-3,\ldots\}\cup\{-\beta,-\beta -1,-\beta-2,\ldots\}\right)$},\\
  0 & \mbox{if $x \in \{-1,-2,-3,\ldots\}$}.
 \end{array} \right.
\label{weight-2-kind}
\end{equation}
The monic multiple Meixner polynomials of second kind, corresponding to the multi-index $\vec{n}$ and the parameters $\vec{\beta}=(\beta_1,\ldots,\beta_r)$, $\beta_i>0$ ($\beta_i-\beta_j \notin \mathbb{Z}$ for all $i\not=j$) and $0<c<1$, is the unique monic polynomial
$M_{\vec{n}}^{\vec{\beta};c}$ of degree $|\vec{n}|$ that satisfies the orthogonality conditions
\begin{equation*}
\sum_{k=0}^{+\infty}M_{\vec{n}}^{\vec{\beta};c}(k)(-k)_jw^{\beta_i}(k)=0,
\qquad j=0,\ldots,n_i-1,\quad i=1,2,\dots,r.
\end{equation*}

For this multiple orthogonal polynomial the Rodrigues formula is (see \cite{arvesu_vanAssche})
\begin{equation*}  \label{MIIRod}
M_{\vec{n}}^{\vec{\beta};c}(x)=
\left(\frac{c}{c-1}\right)^{|\vec{n}|}
\prod_{k=1}^{r}(\beta_i)_{n_i}
\frac{\Gamma(x+1)}{c^x}\prod_{i=1}^{r}\left(\frac{\Gamma(\beta_i)}{\Gamma(\beta_i+x)}\nabla^{n_i}
\frac{\Gamma(\beta_i+n_i+x)}{\Gamma(\beta_i+n_i)}\right)
 \frac{c^x}{\Gamma(x+1)}.
\end{equation*}

In particular, for $r=2$ one gets
\begin{eqnarray*}
\kern-10pt
M_{n_1,n_2}^{\beta_1,\beta_2;c}(x) & = & \left(\frac{c}{c-1}\right)^{n_1+n_2}(\beta_2)_{n_2}(\beta_1)_{n_1}
F_{1:0;1}^{1:1;2}
\left(
\begin{array}{cc}
\begin{array}{c}
(-x):(-n_1);(-n_2,\beta_1+n_1);\\[2ex]
(\beta_1):-;(\beta_2);
\end{array}
\frac{c-1}{c},\frac{c-1}{c}
\end{array}
\right)\\
& = &
\left(\frac{c}{c-1}\right)^{n_1+n_2}(\beta_2)_{n_2}(\beta_1)_{n_1}
\sum_{j=0}^{n_1+n_2}\sum_{k=0}^j
\frac{(-n_1)_k(-n_2)_{j-k}(\beta_1+n_1)_{j-k}}{k!(j-k)!(\beta_2)_{j-k}}
\frac{\left(\frac{c-1}{c}\right)^{j}}{(\beta_1)_{j}}(-x)_{j},
\end{eqnarray*}
where
\begin{eqnarray*}
F_{l:m;n}^{p:q;k}
\left(
\begin{array}{cc}
\begin{array}{c}
\vec{a}:\vec{b};\vec{c};\\[2ex]
\vec{\alpha}:\vec{\beta};\vec{\gamma};
\end{array}
x,y
\end{array}
\right)
& = & \sum_{r,s=0}^{+\infty}\frac{\prod\limits_{j=1}^p(a_j)_{r+s}\prod\limits_{j=1}^q(b_j)_r\prod\limits_{j=1}^k(c_j)_s}
{\prod\limits_{j=1}^l(\alpha_j)_{r+s}\prod\limits_{j=1}^m(\beta_j)_r\prod\limits_{j=1}^n(\gamma_j)_s}
\frac{x^r}{r!}\frac{y^s}{s!},
\end{eqnarray*}
with $\vec{a}=(a_1,\ldots,a_p)$, $\vec{b}=(b_1,\ldots,b_q)$,
$\vec{c}=(c_1,\ldots,c_k)$, $\vec{\alpha}=(\alpha_1,\ldots,\alpha_l)$,
$\vec{\beta}=(\beta_1,\ldots,\beta_m)$,
$\vec{\gamma}=(\gamma_1,\ldots,\gamma_n)$ are the generalizations of the
Kamp\'e de F\'eriet's series \cite{Srivastava} which are a generalization of
the four Appell series in two variables. From this explicit expression can be
found  the coefficients of the recurrence relation.

\begin{theorem}\label{1.2}{\bf(see \cite{arvesu_vanAssche})}
Multiple Meixner polynomials (second kind) $M_{n_1,n_2}^{\beta_1,\beta_2;\,c}$
satisfy the following recurrence relations
\[xM^{\beta_1,\beta_2;\,c}_{n_1,n_2}(x)=M^{\beta_1,\beta_2;\,c}_{n_1+1,n_2}(x)+
b_{n_1,n_2}M^{\beta_1,\beta_2;\,c}_{n_1,n_2}(x)+c_{n_1,n_2}M^{\beta_1,\beta_2;\,c}_{n_1,n_2-1}(x)+
d_{n_1,n_2}M^{\beta_1,\beta_2;\,c}_{n_1-1,n_2-1}(x),
\]
where coefficients are
\begin{eqnarray}
b_{n_1,n_2}  & = &
n_1(2a+1)+n_2(a+1)+a\beta_1, \nonumber \\
c_{n_1,n_2} & = &
a(a+1)\Bigl(n_1n_2+n_1(n_1+\beta_1-1)+n_2(n_2+\beta_2-1)\Bigr),
  \label{MIIrec} \\
d_{n_1,n_2} &
= &  a^2(a+1)n_1(n_1+\beta_1-1)(n_1+\beta_1-\beta_2),  \nonumber
\end{eqnarray}
with $a=\frac{c}{1-c}$.
\end{theorem}

\section{Equilibrium relations and algebraic functions}\label{equil-alg}

Inspired in \cite{RaDiscr1} and \cite{SD} we introduce here the main notions
of the potential theory approach to the weak asymptotics (starting from the
discrete measure of orthogonality) of the polynomials  $ P_n(x), $ orthogonal
with respect to a discrete orthogonality measure \eqref{discrmes}. Usually one
starts with scaling of the problem, i.e.
$$x \to t:\,\,x=x(t,n),$$ (for example, $x=t n$) such that the zeros of the
polynomial sequence $ \tilde{P}(t):=P_n(x(t,n)), $ remain in a compact set
$t\in K \Subset \mathbb{C}$ while $n \to \infty$. This scaling produces a
dependence of the orthogonality weight on the parameter $n$ (\textit{varying
weight}):
$$
\rho_n(t)=\rho(x(t,n))=:\exp\{-nV(t)+O(1)\}.
$$
Accordingly, a dependence on $n$ also occurs for its discrete support
$$\{t_{k,n}\}\,:\, x_k=x_k(t_{k,n},\,n),~k\in\mathbb{N},$$ for which one can
define a counting measure for mass points of the support
$$
\tau_n(t)\,:=\,\frac1n\sum_k\,\delta(t\!-\!t_{k,n}),
$$
and the weak limit of this measure
$$
\tau(\,{\!t\!}\,)\,=\,\lim_n \tau_n(t)\,.
$$
Under some technical conditions (for details see \cite{RaDiscr1} and \cite{SD})
the limiting measure
$$
\nu_{P_n}(x(t,n))\to\lambda(t),
$$
minimizes the following energy functional (with external field $V$))
\begin{equation}\label{ExtrEn}
\Ee(\mu):=\int\left(\Pp^\mu(t)+\frac12 V(t)\!\right)d\mu(t),
\end{equation}
among measures from the class
\begin{equation}\label{Constr}
\mu\,:\quad |\mu|=1\,,~\mu\,\leqslant\,\tau\,.
\end{equation}
The condition of boundedness of the measure $\mu$ by measure $\tau$ reflects
an interesting peculiarity of the discrete orthogonality -due to the interlacing of masses of the discrete measure and zeros of the polynomials.

The potential of the extremal measure $\lambda$ satisfies the equilibrium relations
with the presence of the external field
\begin{equation}\label{EquilCon} \Pp^\lambda(t)
\,+\,V(t)\,
\begin{cases}
~\leqslant\kappa, &t\in\supp\lambda,\\
~=\!\kappa, &t\in\supp(\tau-\lambda)\cap\supp(\lambda):=\Sigma,\\
~\g\kappa, &t\in\supp(\tau-\lambda).
\end{cases}\end{equation}
Set $\Sigma$, where the combination of the potential and the external field, is
equal to a constant is called \textit{equilibrium zone}. A set where the
equilibrium measure $\lambda$ touches its constrain, i.e. measure $\tau$, is
called \textit{saturation zone}. As it follows from  equilibrium relations
\eqref{EquilCon} in the saturation zone the combination of the potential and
the external field is strictly less than the equilibrium constant $\kappa$.

To describe the weak asymptotics of multiple orthogonal polynomials an important
notion of \textit{vector potential equilibrium} was introduced by Gonchar and
Rakhmanov in \cite{17} and \cite{18}.

Let $\vec{\Delta}:=\{\Delta_1,\ldots,\Delta_r\}$ be a collection of compact
sets in $\mathbb{C}$ and let
\begin{equation*}
\label{InterMatrix}
D=(d_{i,j})_{i,j=1}^r,
\end{equation*} be a real symmetric nonsingular positive
definite matrix. An additional condition on $D$ to be compatible with
$(\Delta_1,\ldots,\Delta_r)$ is that $d_{i,j} \geq 0$ whenever $\Delta_i \cap
\Delta_j \neq \emptyset$.

For a vector of measures
\[  \vec{\mu} = (\mu_1,\ldots,\mu_r), \quad \supp \mu_j \subset \Delta_j,\quad j=1,\ldots,r\,, \]
the energy functional $I(\vec{\mu})$ is defined as
\begin{equation}  \label{eq:1.12}
  I(\vec{\mu}) = \sum_{i=1}^r\sum_{j=1}^r d_{i,j} I(\mu_i,\mu_j),
\end{equation}
where $I(\mu_i,\mu_j)$ is the mutual energy of two scalar measures
\[    I(\mu_i,\mu_j) = \int_{\Delta_i} \int_{\Delta_j} \ln \frac{1}{|x-t|} \,
    d\mu_i(x)\,d\mu_j(t) . \]
The extremal vector measure $\vec{\lambda}$, minimizing the energy functional
(\ref{eq:1.12}) among all $\vec{\mu}$ (with fixed masses of the coordinates)
possesses the equilibrium properties for $j=1,\ldots,r$
\begin{equation}  \label{eq:1.13}
    U_j^{\vec{\lambda}}(x) = \sum_{i=1}^r d_{ij} \Pp^{\lambda_i}(x)
    \begin{cases}
     = \kappa_j, & x \in \supp \lambda_j = \Delta_j^*, \\
     \geq \kappa_j, & x \in \Delta_j \setminus \Delta_j^*.
    \end{cases}
\end{equation}
Here the vector $\vec{U}^{\vec{\lambda}} =
(U_1^{\vec{\lambda}},\ldots,U_r^{\vec{\lambda}})$ is called \textit{vector
potential of the vector valued measure $\vec{\lambda}$ with respect to the
interaction matrix $D$}.

External fields $\vec{V} = (V_1,\ldots,V_r)$ such that $V_j(t)$ is a function on $t\in \Delta_j$, $j=1,\ldots,r$, and constraining for $\vec{\mu}$ measures as follows
\[  \vec{\tau} = (\tau_1,\ldots,\tau_r)\,: \qquad \mu_j\,\leqslant\,\tau_j\,,\quad \supp \,\tau_j \subset \Delta_j,\quad j=1,\ldots,r, \]
can be incorporated in the vector equilibrium problem in order to develop an approach for multiple orthogonality with respect to varying weights or (and)
with respect to discrete weights (in a similar fashion that \cite{RaDiscr1} and
\cite{SD}). Thus, the energy functional \eqref{ExtrEn} of vector measure with
external fields becomes (see \eqref{eq:1.12})
\begin{equation*}\label{VecEnergExF}
\Ee(\vec{\mu})\,:=\, I(\vec{\mu})\,+\, \frac12 \sum_{j=1}^r\int
V_j(t)\,d\mu_j(t)\,,
\end{equation*}
and extremal vector measure $\vec{\lambda}$, minimizing this energy functional
possesses the equilibrium properties (see \eqref{eq:1.13}) for $j=1,\ldots,r$
\begin{equation}  \label{vecEQcons}
    U_j^{\vec{\lambda}}(x) \,+\, V_j(x)\,\,
    \begin{cases}
    ~\leqslant\kappa_j, &x\in\supp\,\lambda_j,\\
     ~= \kappa_j, & x \in \supp(\tau_j-\lambda_j) \cap \supp \,\lambda_j = \Sigma_j, \\
     ~\geq \kappa_j, & x \in \supp(\tau_j-\lambda_j).
    \end{cases}
\end{equation}

A saddle point extremal problem for functional \eqref{vecEQcons} is considered
in some applications (when positions of the components of $\vec{\Delta}$ are
not known). In such situation one performs the minimization of the energy
functional by measure $\vec{\mu}$ and maximize it by vector $\vec{\Delta}$. As a
result the support of the equilibrium measure  possesses  the following
symmetry relations
  \begin{equation}  \label{Sprop}
   \displaystyle \frac{\partial \left(U_j\,+\,V_j\right)}{\partial n_+} \,= \,
   \frac{\partial \left(U_j\,+\,V_j\right)}{\partial n_-},\quad
         \textrm{on $\Sigma_j$},\,\,j=1,\ldots,r\,,
  \end{equation}
  where $\partial/\partial n_{\pm}$ denotes the normal derivatives on the respective contours. This property is called S-symmetry after Stahl, who introduced it in \cite{St1} and \cite{St2}.

To conclude this short survey on the logarithmic potential approach to the weak
asymptotics  we recall one general class of multiple orthogonal polynomials
introduced by  Nikishin in \cite{kn:Nikishin}. Such systems are defined as
follows. Let $\sigma_1, \sigma_2$ be two finite Borel measures with constant sign, whose supports $\supp{\sigma_1},$ $ \supp{\sigma_2}$ are contained in non
intersecting intervals $\Delta_1,\Delta_2,$ respectively, of the real line
${\mathbb{R}}$. Set
\[ d \langle \sigma_1, \sigma_2 \rangle(x) = \int \frac{d
\sigma_2(t)}{x-t}\, d \sigma_1(x)\,.
\]
This expression defines a new measure  whose support coincides with that of
$\sigma_1$. In general, let $\vec{\sigma} = (\sigma_1,\ldots,\sigma_r)$ be a system of finite Borel measures on the real line with constant sign and compact support. Let $\Delta_k = [a_k,b_k]$ denote the smallest interval which contains the support of $\sigma_k$. Assume that $\Delta_k\cap\Delta_{k+1}=\emptyset,\;
k=1,\ldots,m-1.$ We say that $\vec{S} = (s_1,\ldots,s_r)=
{\mathcal{N}}(\vec{\sigma})$ is the {\em Nikishin system} generated by
$\vec{\sigma}$, when
\begin{equation} \label{eq:systema}
s_1 = \sigma_1, \quad s_2 = \langle \sigma_1,\sigma_2 \rangle, \ldots , \quad
s_r = \langle \sigma_1, \langle \sigma_2,\ldots,\sigma_r \rangle \rangle.
\end{equation}
In what follows (without loss of generality) we assume $r=2$. So if we
introduce the following notation for the Markov function (Cauchy transform)
\[
\widehat{s}(z):=\int\frac{ds(x)}{z-x}\,\,,
\]
then Nikishin system of two measures is: $s_1 = \sigma_1$, $s_2 = \widehat{\sigma_2} \, \,\,\sigma_1$. These measures have the same support. The Nikishin system is a useful model for
 multiple orthogonal polynomials with respect to a collection of measures
 $(\mu_1,\mu_2)$ with  the same support ($\Delta:=\supp \mu_1=\supp \mu_2 $). For such type of systems an important
 role plays the analytic continuation of the function
\begin{equation}\label{NikW}
\mathcal{F}:=\displaystyle{ \frac{\widehat{\mu_{2}}_{+}-\widehat{\mu_{2}}_{-}}
{\widehat{\mu_{1}}_{+}-\widehat{\mu_{1}}_{-}}
 }\,, \quad \mbox{on}\,\,\Delta\,.
\end{equation}
For the pair $(\mu_1,\mu_2)$ forming a Nikishin system \eqref{eq:systema},
$(\mu_1,\mu_2):=(s_1,s_2)$ we have
$$
\mathcal{F}\,=\,\widehat{\sigma_2}\,\in\,H(\overline{\mathbb{C}} \setminus
\Delta_2).
$$
The limiting zero counting measure for multiple orthogonal polynomial
$P_{\vec{n}},\,\vec{n}:=(n,n) $ with respect to a Nikishin system of measures
$(s_1,s_2)$
$$
\nu_{P_n}\to\lambda\,,
$$
is described by equilibrium problem \eqref{eq:1.13} with the matrix of interaction
\begin{equation}\label{NikMatr}
D:=\left(\begin{array}{cc}
\phantom{-}2 &- 1 \\
-1 &\phantom{-} 2 \\
\end{array}\right),
\end{equation}
and $\lambda=\lambda_1$, where
$|\lambda_1|=2$, $\supp\lambda_1 \in \Delta_1$ and $|\lambda_2|=1$, $\supp\lambda_2\in \Delta_2$. Here the measure $\lambda_2$ is the limiting zero counting measure
$$
\nu_{R_n^{(j)}}\to\lambda_2\,,\qquad j=1,2\,,
$$
for the functions of the second kind
\begin{equation*}\label{2kindFun}
R_n^{(j)}:=\int\frac{P_n(x) ds_j(x)}{z-x}\,,\qquad j=1,2\,.
\end{equation*}

Since the original contribution
\cite{24} by Nikishin there have been published numerous papers devoted to the asymptotics of multiple orthogonal polynomial with respect to a Nikishin system of measures (for most general results and latest surveys see \cite{AptLys} and
\cite{AptKuUMN}).

In this section, starting from the fact that components of the vector of discrete Meixner weights have the same support we assume that these measures inherit properties of the Nikishin system of continuous measures with the same support (which is not developed for the discrete measures yet). This assumption allows us heuristically to state an ad hoc vector equilibrium problem with Nikishin matrix of interaction, with external field (because of unbounded set of orthogonality) and constrain (due to discrete weight). Accordingly, at this stage we look for the solutions of this equilibrium problem in a form of spectral curve (among the simplest candidates satisfying the necessary conditions following from the equilibrium problem). Thus, as a result of the first stage of our analysis we have heuristically got equations of the algebraic curves \eqref{EQphi} and \eqref{phieqII} and have related their logarithms (rigorously)  with certain vector equilibrium problems (see Propositions ~\ref{prop:2.1} and \ref{prop:2.2}). We highlight that at this point spectral curves and equilibrium problems are not formally related to the multiple Meixner orthogonal polynomials.

In section~\ref{nroot} we derive the same spectral curves from the recurrence coefficients of the multiple Meixner orthogonal polynomials. Consequently, it now rigorously relates Nikishin type equilibrium with constrains to Meixner discrete orthogonality.

\subsection{Potential problem and its solution for classical Meixner polynomials }\label{Ss2.1}

It this subsection, aimed to illustrate our approach in connection with the above discussed techniques, we take (as example) the classical Meixner polynomials~\eqref{meixpol}. First, we state the equilibrium problem for the limiting zero counting measure. Then, we show our approach, which is based on an algebraic function formulation, by solving this equilibrium problem in a quite direct way.

The scaling of variable $x$ in \eqref{meixmes}
$$
x=:nt,
$$
stabilizes the zeros of the polynomials (i.e. when $n\!\to\!\infty$, the zeros stay in a compact)
$$
\~P_n(t):=M_n(x),
$$
which one considers now as a polynomial of variable $t$. Indeed, it is orthogonal to $\{t^k\}_{k=0}^{n-1}$ with respect to the varying weight (depending on $n$)
$$
\~\rho_n(t):=\frac{\Gamma(\beta+tn)}{(tn)!}\,c^{tn}\;,
$$
concentrated in the points $\{\frac kn\}_{k\in\mathbb{N}}$. Taking the limit
$$
\lim_{n\to\infty}|\~{\rho}_n(t)|^{1/n}=c^{t}=\exp\left\{-t\,\ln\frac1c\right\}\;,
$$
we see that the external field for extremal problem
\eqref{ExtrEn}\,--\,\eqref{Constr} and equilibrium problem \eqref{EquilCon} has
the follownig form
\begin{equation}\label{1.7}
V(t)=\ln\frac1c\;\R\,t\;.
\end{equation}
We also have
\begin{equation}\label{1.8}
\frac1n\,\sum\limits_{k=0}^\infty\delta\left(t-\frac kn\right) \stackrel
*\longrightarrow\,l(t)\,=:\tau(t),\qquad t\in\RR_+\;,
\end{equation}
where $l(t)$ is the Lebesgue measure with unit density $(dt)$ on the semiaxis
$\mathbb{R}_+$.

Thus, the limiting measure $\lambda(t)$ of zero distributions of polynomials
 $\~P_n(t)$
\begin{equation*}\label{1.85}
\nu_{\~P_n(t)}\stackrel *\longrightarrow\lambda(t),
\end{equation*}
satisfy the extremal problem  \eqref{ExtrEn}\,--\,\eqref{Constr} as well as the conditions of equilibrium \eqref{EquilCon} with external field \eqref{1.7} and constrain \eqref{1.8}.

Let us show our approach for finding the explicit representation of
the measure $\lambda(t)$. We look for the Cauchy transform of measure $\lambda$ (or Markov function)
\begin{equation}\label{1.9}
h(z):=\hat{\lambda}(z)=\int\frac{d\lambda(t)}{z-t}\;,\qquad
z\in\mathbb{C}\,\backslash\supp\lambda\;,
\end{equation}
as function $h(z)$, which be analytic on some two sheeted Riemann
surface $\Rr:=\left(\mathbb{C}^{(0)},\mathbb{C}^{(1)}\right)$. The possibility
of analytic continuation of the function $h(z)$ to another sheet of the Riemann
surface through the compact $\Sigma$ can be seen from the second relation in
\eqref{EquilCon} after it ``complexification" and differentiation. Thus,
the function $h:=(h_0,h_1)$ on $\Rr$ has some natural condition for its values at
 $\infty$.
For the main branch $h_0$, because of \eqref{1.9} (and normalization of measure
$|\lambda|=1$)
\begin{equation}\label{1.10}
h(z):=\frac{1}{z}+\cdots\;,\quad z\!\to\!\infty^{(0)}\;,
\end{equation}
and for its continuation $h_1$ (taking into account derivative of the
``complexification" of external field \eqref{1.7})
\begin{equation}\label{1.11}
h(z):=\ln\frac1c-\frac{1}{z}+\,\cdots\;,\quad z\!\to\!\infty^{(1)}\;.
\end{equation}
The density of the measure $\lambda'(t)$ is proportional to the jump of imaginary part of the function $h(z)$ along $\RR_+$. Moreover, it is clear from electrostatic consideration of the problem \eqref{ExtrEn}\,--\,\eqref{Constr} that at the neighborhood of the left end point of $\supp\lambda$, i.e. to the right from the origin, this jump has to be a constant -here the measure  $\lambda$ reaches its constrain by the Lebesgue measure $l(t)$. Such behavior can be modelled by representing the function $h(z)$ as logarithm of some other function, say $E(z)$, which on this part of the support of the measure $\lambda$ (i.e. in the saturation zone) takes negative real values. Hence, we are looking  for $h(z)$ in the form
\begin{equation}\label{1.12}
h(z)=\ln E(z),
\end{equation}
where $E(z)$ is a rational function on $\Rr$ with one zero and
one pole (at the points on the different sheets of Riemann surface $\Rr$ with
the same projection point~0)
\begin{equation}\label{1.13}
E(z)=\left\{\begin{array}{ll} 0\,, & z=0^{(0)},\\
\infty\,, & z=0^{(1)}.
\end{array}\right.
\end{equation}
We have chosen the position of these points in order to provide holomorphicity
of branches $h$ -equation \eqref{1.12} in $\CC\backslash\RR_+$. The conditions
\eqref{1.10} and \eqref{1.11} give the following behavior for $E(z)$ in the
neighborhood of infinity on both sheets of $\Rr$
\begin{equation}\label{1.14}
E(z)=\left\{\begin{array}{ll} 1+\dfrac{1}{z}\,, &z\!\to\!\infty^{(0)},\\
[6pt] \dfrac1c\left(1-\dfrac{1}{z}\right)\,, & z\!\to\!\infty^{(1)}.
\end{array}\right.
\end{equation}
Observe that it is enough to have conditions \eqref{1.13} and \eqref{1.14} for
determining $E(z)$ and equation of its Riemann surface. Indeed, the Vieta
relations give
\begin{equation*}\label{1.151}
E^2+\left[\frac{1}{z}\Big(\frac1c-1\Big)-
\Big(\frac1c+1\Big)\right]\!E+\frac1c=0.
\end{equation*}
From here it follows that Riemann surface of the function $E(z)$ has branch
points
$$
e_1=\frac{1-\sqrt{c}}{1+\sqrt{c}}\;,\quad e_2=\frac1{e_1}\;.
$$
Moreover, both branches of the function $E(z)$ are negative on $[0,e_1]$ and
positive on  $[e_2,\infty]$. Therefore, the jump of the imaginary part of the function $h(z)$ -see equation \eqref{1.12}- is a constant on $[0,e_1]$ and equal zero on  $[e_2,\infty]$. Thus, $\mbox{supp}\,\lambda=[0,e_2]$, the interval $[0,e_1]$ is saturation zone, the interval $[e_1,e_2]=\Sigma$ corresponds to the equilibrium zone, and $[e_2,\infty]$ is a zone free from zeros of the scaled classical Meixner polynomials. Measure $\lambda$ obtained by this procedure satisfies the equilibrium condition \eqref{EquilCon}. Furthermore, after integration one can obtain an explicit value of the equilibrium constant in \eqref{EquilCon}
$$
\kappa=1+\ln\frac{1-c}c.
$$

\subsection{Potential problem and its solution for multiple Meixner polynomials (first kind)}
\label{Ss2.2}

Assume $r=2$, and $\vec{n}=(n,n)$. We rewrite  weights
(\ref{weights-kind-1}) as follows
$$
w_{j}^{\beta}(x)=\frac{\Gamma(\beta+x)}{x!}e^{x\ln
c_{j}},\quad  j=1,2.
$$
Here the supports of both discrete orthogonality measures coincide. Therefore, we pose an equilibrium problem for the weak asymptotics
in a similar way to the Nikishin case. Indeed, we consider the equilibrium
problem associated with the weights $w_{1}^{\beta}(x)$, and
$u(x)w_{1}^{\beta}(x)$, where
$$
u(x)=\frac{w_{2}^{\beta}(x)}{w_{1}^{\beta}(x)}=e^{x\ln
\frac{c_2}{c_1}} =\frac{1}{2\pi
i}\int_{\Gamma}\frac{e^{t\ln\frac{c_2}{c_1}}}{t-x}dt.
$$
The integral representation in the right hand side is due to the Cauchy integral
formula with a contour of integration $\Gamma$ encircling $\mathbb{R}_+$, i.e.
it comes from infinity through the upper half plane along the positive semiaxis
turning counterclockwise around zero in the left half plane and returning to infinity along the positive semiaxis through the lower half plane. In addition, the contour $\Gamma$ serves as a boundary of holomorphicity region of an analog of the function \eqref{NikW}.

Now we do a scaling of the problem in order to compactify zeros of   multiple
Meixner polynomials (first kind). Let
$$\widetilde{M}_{2n}^{\beta;\vec{c}}(x)=C_{n}M_{n,n}^{\beta;\vec{c}}(xn)=x^{2n}+\cdots\,.$$
Recall that the function of the second kind associated with the multiple
Meixner polynomials of first kind is given by
$$
R_n^{(j)}(z)=\int\frac{M_{n,n}^{\beta;\vec{c}}(t)}{t-z}d\mu_j,\quad j=1,2\,\,,
$$
where $\mu_j$ is the discrete orthogonality measure \eqref{mesMeixI}. Let us
denote $\widetilde{R}_{n}^{(j)}(z)=R_n^{(j)}(nz)$.
General feature of the Nikishin systems is that the functions of the second kind have $n$-zeros on a compact of the complex plane outside of the support of the orthogonality measures. For the classical Nikishin system (see \eqref{eq:systema}) they accumulate along interval $\Delta_2$. For the case under consideration we expect that the function
$\widetilde{R}_n^{(j)}(z)$ has $n$-zeros accumulated along an analytic arc
$\gamma$ in the complex plane. We assume that contour $\Gamma$ contains this
analytic arc~$\gamma$.

For the formulation of the equilibrium problem, we
follow the Nikishin system model, in which the zeros accumulate along an open
contour $\gamma\subset\mathbb{C}$ (for the functions of the second kind) and in
$\mathbb{R}_+$ (for the multiple orthogonal polynomials). Since the mass points
$x_k$ are uniformly distributed over a subset of $\mathbb{R}$, then the
zero counting measure
$\nu_{\widetilde{M}_{2n}^{\beta;\vec{c}}(x)}\overset{*}{\to}\lambda$, where
$\abs{\lambda}=2$, and $\lambda\leq l$, being $dl(x)=dx$ the Lebesgue measure,
and $\supp(\lambda)\subset\mathbb{R}_{+}$. On the other hand,
$\nu_{\widetilde{R}_{n}^{(j)}(x)}\overset{*}{\to}\mu$, and its magnitude is
equal to $1$ with $\supp(\mu)\subset\gamma$. Accordingly,
\begin{equation}\label{VmesMMI}
 \begin{cases}
\,\,\lambda\,: \quad \abs{\lambda}=2,
\quad \lambda\leq l,\qquad & \supp(\lambda)\subset\mathbb{R}_{+},\\
\,\,\mu\,:\quad \abs{\mu}=1,\quad \quad & \supp(\mu)=:\gamma\subset\Gamma.
\end{cases}
\end{equation}
Now  we write the equilibrium problem for these measures
$(\lambda,\mu):=\vec{\lambda}$ as \eqref{vecEQcons}, where:
\begin{itemize}
 \item The matrix of interaction $D$ in the definition \eqref{eq:1.13} of the
 vector potential $U_j^{\vec{\lambda}}$ has the Nikishin's form \eqref{NikMatr}.

\item Scaling of the weights $w_1^\beta$ and $u$ brings external fields
$V_1:=(-\ln c_1)\R{x}$ on $\mathbb{R}_+$ and $V_2:=\left(-\ln
\frac{c_2}{c_1}\right)\R{x}$ on $\Gamma$, respectively.
\end{itemize}
Thus, we have the equilibrium relations:
\begin{equation}\label{EqConMMI}
 \begin{cases}
\,\,W_1\,:=\,2\mathcal{P}^{\lambda}-\mathcal{P}^{\mu}-(\ln
c_1)\R{x}&\left\{
\begin{array}{l}\leq\kappa_{1},\quad\mbox{on}\quad\supp(\lambda)\subset\mathbb{R}_{+},\\
\geq\kappa_{1},\quad\mbox{on}\quad\mathbb{R}_{+}\cap\supp(l\setminus\lambda).
                                        \end{array}\right.
\\
\,\,W_2\,:=\,2\mathcal{P}^{\mu}-\mathcal{P}^{\lambda}-\left(\ln
\frac{c_2}{c_1}\right)\R{x}&\left\{
\begin{array}{l}=\kappa_{2},\quad\mbox{on}\quad\supp(\mu)=:\gamma\subset\Gamma,\\
                                          \geq\kappa_{2},\quad\mbox{on}\quad\Gamma.
                                      \end{array}\right.
\end{cases}
\end{equation}
We also have to add to this system the $S$-symmetry property \eqref{Sprop}  for
the second plate $\Gamma$ of our condenser $(\mathbb{R}_+, \Gamma)$:
\begin{equation}  \label{SpropMI}
   \displaystyle \frac{\partial W_2}{\partial n_+} \,= \,
   \frac{\partial W_2}{\partial n_-},\quad
         \textrm{on $\gamma$}.
  \end{equation}

Now we look for a solution of the equilibrium problem \eqref{VmesMMI},
\eqref{EqConMMI}, \eqref{SpropMI}. In accordance with our approach (by
differentiating) we pass from
 \eqref{EqConMMI} to the following relations for the  real part of the Cauchy transforms of the measures $\lambda$ and $\mu$,
\begin{align*}
\R\left\{2\hat{\lambda}-\hat{\mu}+\ln c_1\right\}&=0,\quad
\mbox{on}\quad \supp(\lambda)\cap\supp(l\setminus\lambda)=:\Sigma,\\
\R\left\{2\hat{\mu}-\hat{\lambda}+\ln c_2-\ln c_1\right\}&=0,\quad\mbox{on}
\quad \supp(\mu):=\gamma,
\end{align*}
respectively. From these expressions one rewrites
\begin{align}\label{conmMI}
\hat{\lambda}&=\hat{\mu}-\hat{\lambda}-\ln c_1,\quad \mbox{on}\quad \Sigma\,,\\
\hat{\mu}+\ln c_2&=\hat{\lambda}-\hat{\mu}+\ln c_1,\quad \mbox{on}\quad
\gamma\,.\notag
\end{align}
Notice that if we define on the three sheeted Riemann surface
$$
\mathcal{R} \,:=\, \left(\mathcal{R}_0, \mathcal{R}_1, \mathcal{R}_2
\right)\,:\qquad
\begin{cases}
\,\,\mathcal{R}_0\,:=\,\overline{\mathbb{C}}\setminus \Sigma\,,\\
\,\,\mathcal{R}_1\,:=\,\overline{\mathbb{C}}\setminus \left\{\,\Sigma\bigcup \gamma\right\}\,,\\
\,\,\mathcal{R}_2\,:=\,\overline{\mathbb{C}}\setminus  \gamma,
\end{cases}
$$
a function
$$
H(z):=\left\{\begin{array}{ll}
\hat{\lambda},\quad & \mbox{in}\quad \mathcal{R}_0,\\
\hat{\mu}-\hat{\lambda}-\ln c_1,\quad & \mbox{in}\quad \mathcal{R}_{1},\\
-\hat{\mu}-\ln c_2,\quad & \mbox{in}\quad
\mathcal{R}_{2},\end{array}\right.
$$
then due to \eqref{conmMI} and \eqref{SpropMI} the branches of this function
have analytic continuation one to another through the cuts joining the sheets
$\left(\mathcal{R}_0, \mathcal{R}_1, \mathcal{R}_2 \right)$ of the Riemann
surface $\mathcal{R}$. Observe that from \eqref{VmesMMI} follows that $H(z)$ behaves near infinity as
$$
H(z)=\left\{\begin{array}{ll}
2z^{-1}+\cdots,\quad&\mbox{as $z\to\infty\,\,$ in}\quad \mathcal{R}_{0},\\
-\ln c_1-z^{-1}+\cdots,\quad&\mbox{as $z\to\infty\,\,$ in}\quad \mathcal{R}_{1},\\
-\ln c_2-z^{-1}+\cdots,\quad&\mbox{as $z\to\infty\,\,$ in}\quad
\mathcal{R}_{2}.\end{array}\right.
$$
In the same way like we did it for the classical Meixner polynomials we model
the constrain condition  $\lambda\leq l$ and $dl(x)=dx$ on $\mathbb{R_+}$,
representing the function $H(z)$ as logarithm of a meromorphic function $\phi(z)$ on $\mathcal{R}$, i.e.
\begin{equation}\label{HnMI}
H(z)=\ln\phi(z)\,,\qquad \phi \in \mathcal{M}(\mathcal{R}).
\end{equation}
To have a jump for the function $H$ along $\mathbb{R_+}$ lifted to sheet
$\mathcal{R}_{0}$ we force the function $\phi(z)$ to take values $\phi=\infty$, and $\phi=0$ at the points  from the two sheets  $\left(\mathcal{R}_0,
\mathcal{R}_1\right)$ which have on $\mathbb{C}$ projection point $z=0$
\begin{equation}
\phi(z)=\left\{\begin{array}{ll}
\,0\,,\quad&\mbox{as $z\to 0\,\,$ in}\quad \mathcal{R}_{0},\\
\infty,\quad&\mbox{as $z\to 0\,\,$ in}\quad \mathcal{R}_{1}.\end{array}\right.
\label{behavior 0_phi_meixner1}
\end{equation}
Thus, $H \in
\mathcal{M}\left(\mathcal{R}\setminus\left\{\mathbb{R}_+^{(0)}\bigcup\mathbb{R}_+^{(1)}\right\}\right)$.
From the behavior of $H$ at infinity we conclude that
\begin{equation}
\phi(z)=\left\{\begin{array}{ll}
1+2z^{-1}+\cdots,\quad&\mbox{as
$z\to\infty\,\,$ in}\quad
\mathcal{R}_{0},\\
 c_1^{-1}\left(1-z^{-1}\right)+\cdots,\quad&\mbox{as $z\to\infty\,\,$ in}\quad \mathcal{R}_{1},\\
 c_2^{-1}\left(1-z^{-1}\right)+\cdots,\quad&\mbox{as $z\to\infty\,\,$
in}\quad \mathcal{R}_{2}.\end{array}\right. \label{behavior_infty phi_meixner1}
\end{equation}
If $\Sigma$ and $\gamma$ are connected then the genus of the Riemann surface $\mathcal{R}$ is equal to zero (we assume that it is true). The simplest meromorphic function $\phi$ (i.e.
see \eqref{behavior 0_phi_meixner1} -one pole and one zero) that maps
$\mathcal{R}$ to the Riemann sphere and therefore the inverse function $z(\phi)$ is just a rational function. From \eqref{behavior_infty phi_meixner1} we know the poles $z$ as a  function of $\phi$ and moreover we know  its residues in the poles. Taking into account also  \refe{behavior 0_phi_meixner1} one gets that the algebraic function $\phi(z)$ is defined by the equation
\begin{equation}
\phi(z)\,:\qquad
z=\dfrac{2}{\phi-1}-\dfrac{1}{ c_{1}\phi-1}-\dfrac{1}{ c_{2}\phi-1}.
\label{EQphi}
\end{equation}
This algebraic function has four branch points. In order to
construct (for the condenser with two plates $(\mathbb{R}_+,
\Gamma), \,\, \Gamma \subset\mathbb{C}\setminus\mathbb{R}_+$) a
solution of the equilibrium problem \eqref{EqConMMI}-\eqref{SpropMI}  we need to have  the case when:
\begin{itemize}
\item $\phi$ has two complex conjugated branch
points $e, \overline{e}$, and two real $e_1, e_2$ ($0<e_1<e_2$);
\item the analytic arc $\gamma$ joining points $(e,
    \overline{e})$, chosen by means of \eqref{SpropMI}, does
    not cross $\mathbb{R}_+$;
\item values $\phi(\zeta) \in \mathbb{R}_{-}$, when  $\zeta \in \,[0,e_1]^{(0)}\subset \mathcal{R}_{0}$.
\end{itemize}

\begin{figure}[ht!]
\centerline{\includegraphics[width=0.45\textwidth]{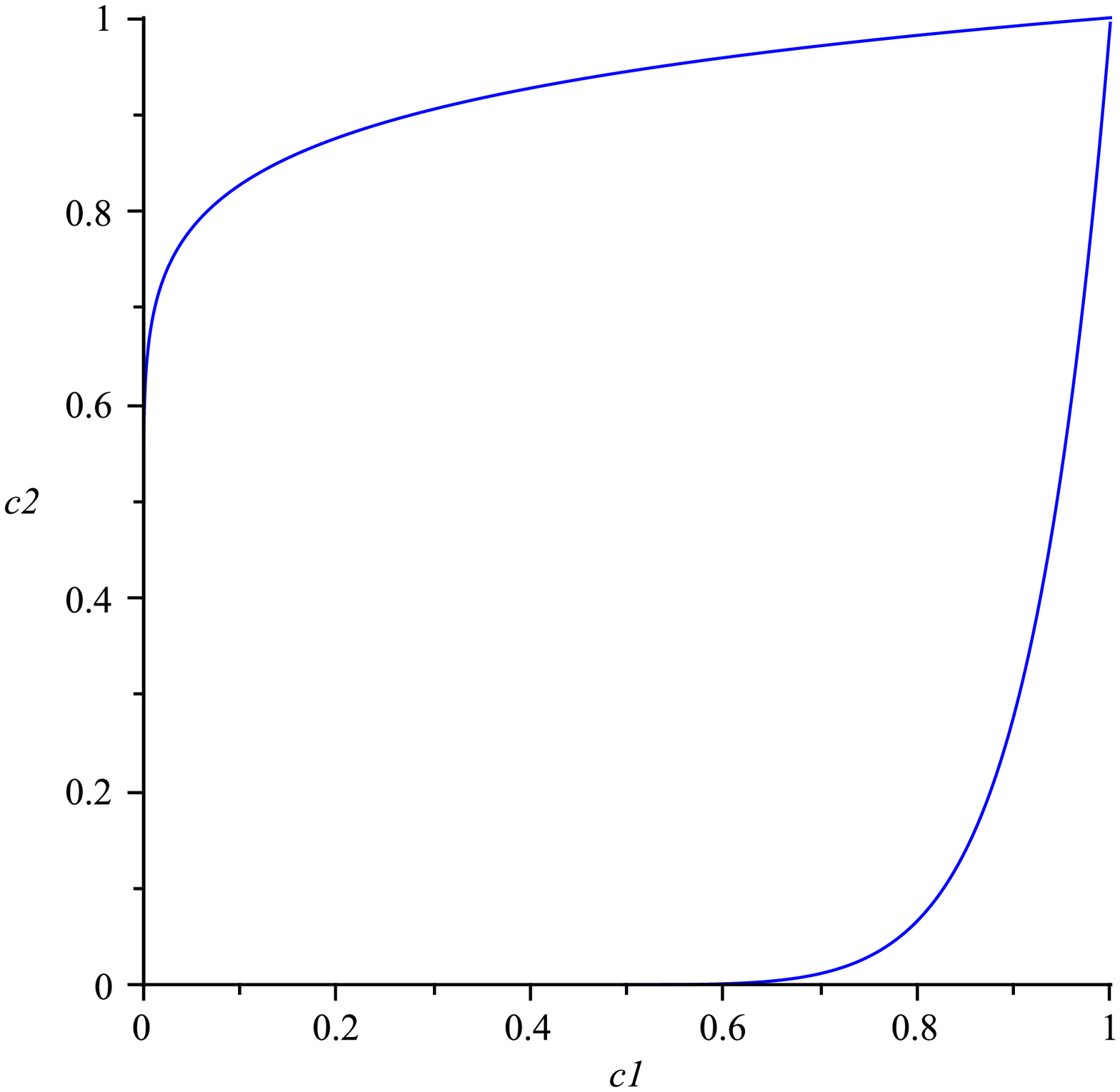}
\qquad
\includegraphics[width=0.45\textwidth]{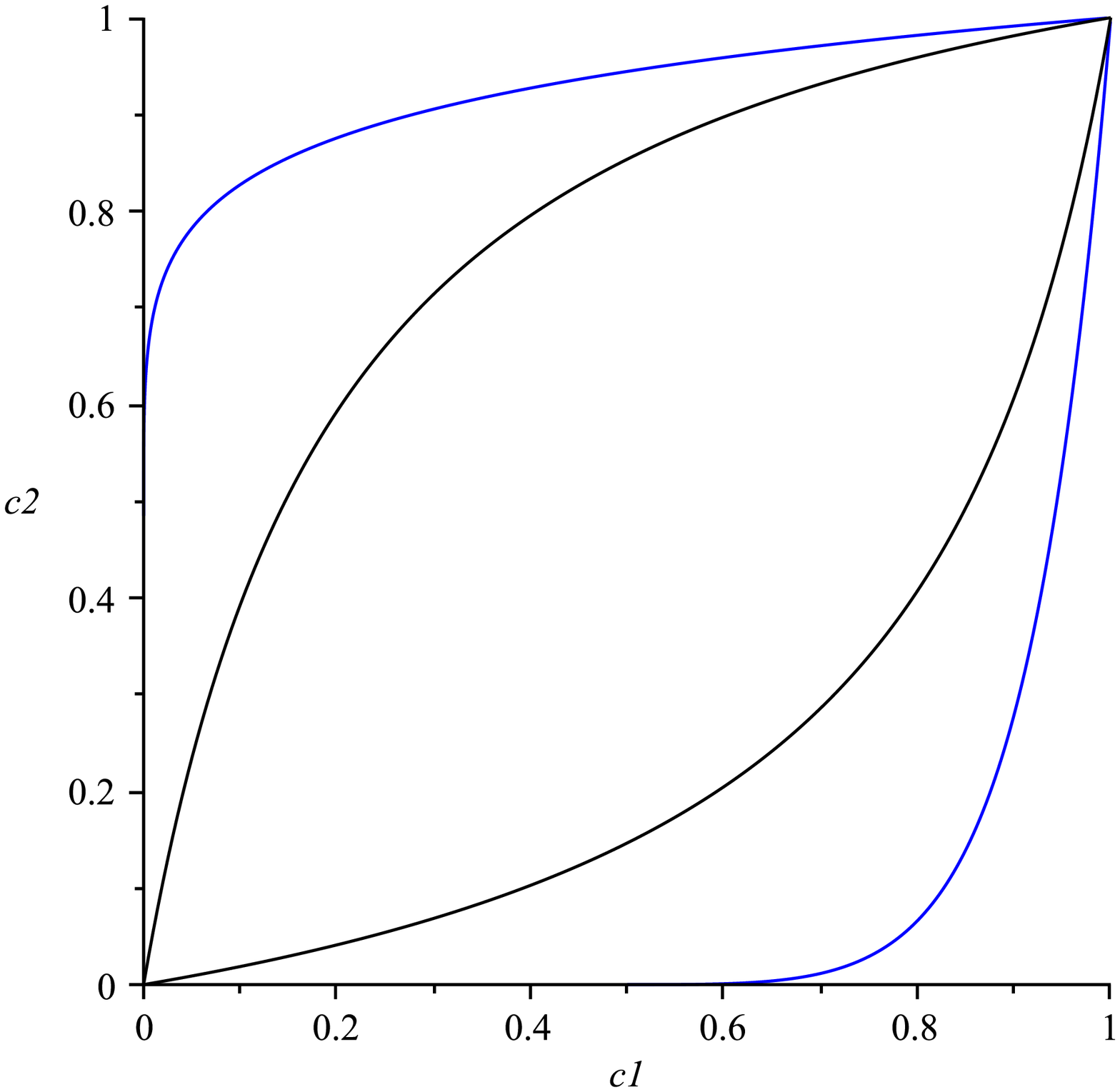}}
\caption{(A) --on the left. Domain $\mathcal{GN}$. If $(c_1,c_2)\in \mathcal{GN}$, then  $\phi$ has two complex  and two real branch
points;
(B) --on the right. Domain $\mathcal{N} \subset \mathcal{GN}$. If $(c_1,c_2)\in \mathcal{N}$
then  $\gamma \cap \mathbb{R}_+ = \varnothing$.} \label{fig1}
\end{figure}

The first item can be checked by analyzing the critical values
$c_{1}, c_{2}$ when the branch points of the function $\phi$ coincide
--the transition of  two conjugate branch points to the case of
two real branch points. For this purpose we take consequently two
times the discriminant of the curve \eqref{EQphi} --first time in
the variable $\phi$ and second time in the variable $z$. In result
we get a polynomial of two variables $c_1$ and $c_2$ with rational
coefficients. Factorization of this polynomial gives a factor of
high degree (maximal degree in each variable $c_{j}$ is 6 and
maximal degree of the monomials is 9). It defines an algebraic
curve $ T(c_{1},c_{2})=0,$ which has genus equal 1. Thus, this curve
is rather complicated. However, taking into account the symmetry of
this curve with respect to $(c_{1},c_{2})$, and changing variables
\begin{equation*}
\left\{\begin{array}{ll}
\,c_{1}+c_{2}&\,=:\,s\,,\\
\,\,\,c_{1}\,\,c_{2}&\,=:\,p\,,\end{array}\right.
\end{equation*}
yields the algebraic curve
\begin{equation}\label{CriticalCurve}
\widetilde{T}(s,p)\,:=\,T(c_1(s,p), c_2(s,p))\,=0\,,
\end{equation}
which has genus equal 0. Therefore, curve $ \widetilde{T}(s,p) $ admits a
rational uniformization. Using Gr\"{o}bner basis technique we arrive to
the following simple parametrization of this curve
\begin{equation*}
\left\{\begin{array}{ll}
\,\,s\,=\,&-\displaystyle\frac{(t^{2}+ 8 t + 18)(t + 2)^{2} }{2t}\,,\\
\,\,p\,=\,&-\displaystyle\frac{(t + 3)^{3}(t + 2)^{3} }{t^{3}}\,.\end{array}\right.
\end{equation*}
An analysis of this curve shows that $c_{1},c_{2}\,\in\,\mathbb{R}$
when $t\,\in\,[-\infty,-2] \cap [0, \infty]$. On Figure~\ref{fig1}~(A)
 for $t\,\in\,[-4.73 ,-3] $ we see branches of the critical curve \eqref{CriticalCurve}
 which belong to the interesting for us quadrate  $(c_1,c_2)\in \{(0,1)\times(0,1)\}$. If parameters $(c_1,c_2)$ are inside of the domain bounded by these branches
 and axes $\{c_1=0, c_2\in(0, 0.5)\}$  and $\{c_2=0, c_1\in(0, 0.5)\}$ then function
 $\phi$ has two complex conjugated branch
points $e, \overline{e}$, and two real $e_1, e_2$ ($0<e_1<e_2$). We denote this domain $\mathcal{GN}$. On the critical lines the conjugate branch points
 $e, \overline{e}$ coincide on the real axes and when  $(c_1,c_2)\in \{(0,1)\times(0,1)\} \setminus \mathcal{GN}=: \mathcal{A}$, the function $\phi$ has four real branch
points.

\begin{figure}[ht!]
\centerline{\includegraphics[width=0.45\textwidth]{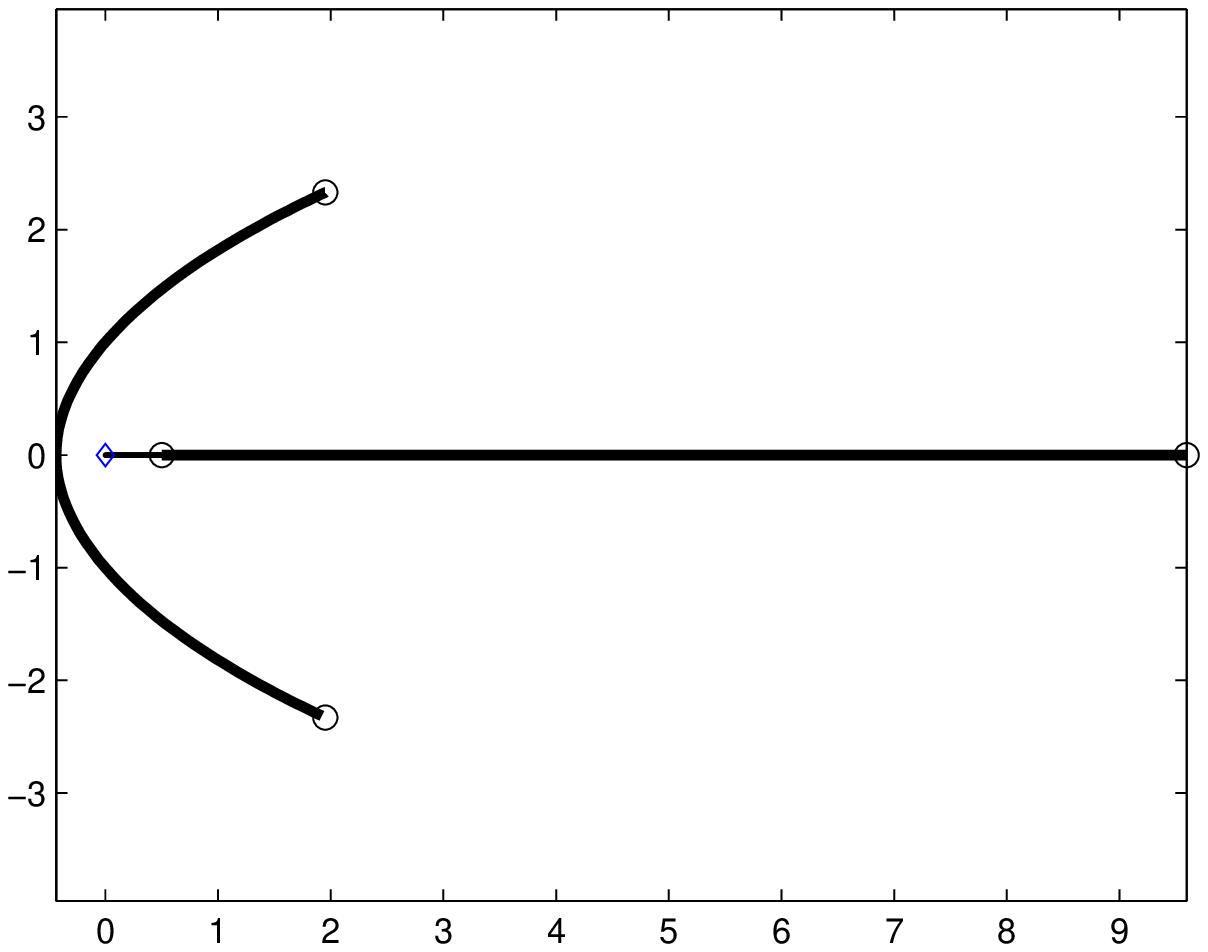}
\qquad
\includegraphics[width=0.45\textwidth]{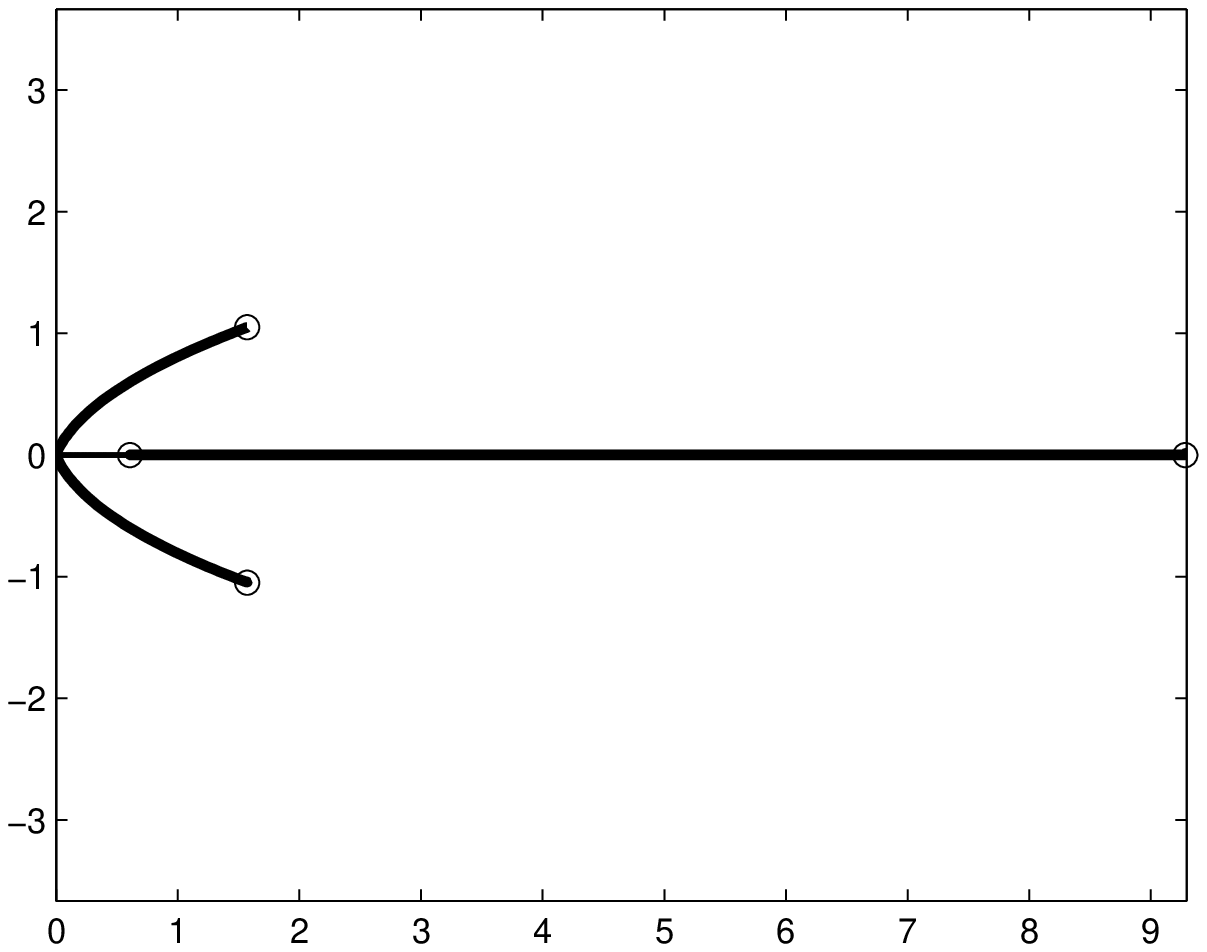}}
\caption{The sets $\gamma = \mbox{supp}(\mu)$ and $[0,e_2] = \mbox{supp}(\lambda)$
as well as the branch points (circles) are depicted for the cases:
(A) --on the left. $(c_1,c_2)=(0.5, 0.25)\in \mathcal{N};\quad$
(B) --on the right. $(c_1,c_2)=(0.5, 0.146445649)$ is critical.}\label{fig2}
\end{figure}

The second item, i.e. condition on the parameters $(c_1,c_2)$ implying  $\gamma \cap \mathbb{R}_+ = \varnothing$ can be clarified using an equivalent (constructive) definition of the  S-curve $\gamma$ (see \eqref{SpropMI}). It is known (see \cite{Lys}, \cite{St1}) that S-property condition \eqref{SpropMI} means that the analytic curve $\gamma$ belongs to the union of the critical (orthogonal) trajectories of the quadratic differentials
\begin{equation}\label{QD}
\R \int (H_j-H_k)\, d z \,=\,0\quad\Longleftrightarrow\quad (H_j-H_k)^2 (d z)^2\,<\,0.
\end{equation}
Using uniformization \eqref{HnMI}, \eqref{EQphi} we perfom integration by parts
$$ \int H(z)\, d z \,=\,z\,\ln \phi\,-\,\int \frac{z}{\phi}\,d\phi\,=\,z\,\ln \phi\,+\,\frac{(\phi-1)^2}{(\phi-\sigma_1)(\phi-\sigma_2)},\qquad \sigma_j:=1/c_j,\,\,\,j=1,2.$$
Thus condition \eqref{QD} is equivalent to
\begin{equation}\label{QDII}
\left|\exp\left\{\int (H_j-H_k)\, d z \right\}\right|\,=\,1\quad\Longleftrightarrow\quad\left|\left(\displaystyle\frac{\phi_j}{\phi_k}\right)^z\,\displaystyle
\frac{(\phi_j-1)^2(\phi_k-\sigma_1)(\phi_k-\sigma_2)}{(\phi_k-1)^2(\phi_j-\sigma_1)(\phi_j-\sigma_2)}\right|\,=\,1.
\end{equation}
Now we obtain the critical values of the parameters such that curve
$\gamma$ is crossing the origin. Substituting in \eqref{QDII}
 $j=1,\,k=2,\,z=0$ and taking from \eqref{behavior 0_phi_meixner1},
 \eqref{EQphi}  value of $\phi(0)$
 we have
 $$ \left |\displaystyle
\frac{(\phi_2(0)-\sigma_1)(\phi_2(0)-\sigma_2)}{(\phi_2(0)-1)^2}\right|\,=\,1,\qquad
\mbox{and}\qquad \phi_2(0)\,=\,
\frac{2 \sigma_1 \sigma_2 \,
-\,\sigma_1\,-\,\sigma_2}{\sigma_1\,+\,\sigma_2\,-\,2}.
 $$
Therefore, taking into account that $\sigma_1,\sigma_2 \in \mathbb{R}$ we arrive to
$$\left |\displaystyle
\frac{(\sigma_1-1)(\sigma_2-1)}{(\sigma_1-\sigma_2)^2}\right|\,=\,\frac14,
\qquad\mbox{or} \qquad \frac{ (\sigma_1+ \sigma_2-2)^2 \, -\,(\sigma_1\,-\,\sigma_2)^2}{(\sigma_1\,-\,\sigma_2)^2}\,=\,\pm1\,.
$$
A branch of the critical curve corresponding to sign $(-)$ in the
right hand side lies out of the admissible range for $c_{1}, c_{2}$
(or $\sigma_{1}, \sigma_{2}$). Finally, considering in the right
hand side sign $(+)$ we get
\begin{equation}\label{QDcrit}
(\sigma_1\,+\sigma_2-2)^{2}\,=\,2(\sigma_1\,-\sigma_2)^2.
\end{equation}
On Figure~\ref{fig1}~(B) the branches of the critical curve \eqref{QDcrit} forms the boundary of the
domain $\mathcal{N}$. In this domain $(c_1,c_2)\in\mathcal{N}$ the analytic arc
$\gamma$ does not intersect $\mathbb{R}_+$ (see Figure~\ref{fig2}~(A)). On the
boundary of this domain $\gamma$ is passing trough the origin (see Figure~\ref{fig2}~(B)).

 The cut joining the logarithmic singularities of the function
 $H$ starts at  point $0$ on sheet $\mathcal{R}_{0}$, goes along $\mathbb{R}_+$ to branch point $e_1$, passing  to  sheet $\mathcal{R}_{1}$ and goes there
 back to point $0$. Along this cut (for all  $(c_1,c_2)\in \{(0,1)\times(0,1)\}$) values $\phi$ are negative. Thus we have, indeed $$H \in
\mathcal{M}\left(\,\mathcal{R}\,\setminus\,\left\{\,[0,e_1]^{(0)}\bigcup\,[e_1,0]^{(1)}\right\}\right)\,.$$

Summarizing we have that for the parameters belonging to the domain
\begin{equation}\label{NikDom}
\mathcal{N}\,:=\,\left\{(c_1,c_2)\in \{(0,1)\times(0,1)\}\,:\quad \frac{c_2(1-\sqrt{2})}{2 c_2-1-\sqrt{2}}
\,<\,c_1\,<\, \frac{c_2(1+\sqrt{2})}{2 c_2-1+\sqrt{2}}\,\right\},
\end{equation}
the density of  measure $\lambda$ is determined by the limiting
values at the points $[0,e_2]$ of the imaginary part of  $\ln\phi(z)/\pi$
 or equivalently,
\begin{equation}\label{EqMesMMI}
\frac{d\lambda(x)}{dx}\,=\,\lim_{y\to 0}\,\frac{1}{\pi}\, \Ima \,H(x+i
y)\,,\qquad x\in [0,\, e_2]\quad \mbox{in}\quad \mathcal{R}_{0}\,.
\end{equation}
The interval $[0,\, e_2]$ is the support of measure $\lambda$.
Since the function $\phi$ is negative on $[0,\, e_1]$ in $\mathcal{R}_{0}$ we have that the density of the measure $\lambda$ is equal to unity
there. Therefore, the interval $[0,\, e_1]$ corresponds to the saturation zone and the interval $\Sigma= [e_1,\, e_2]$ is the equilibrium zone.
The support of measure $\mu$ is the analytic arc $\gamma$.
Thus, we have proved the following result.
\begin{proposition} \label{prop:2.1}
The density of the absolute continues measure $\lambda$,  solving equilibrium
problem \eqref{VmesMMI}, \eqref{EqConMMI}, \eqref{SpropMI} for $(c_1,c_2)\in \mathcal{N}$, is equal to
\eqref{EqMesMMI}, where the function $H$ is defined by \eqref{HnMI} and
\eqref{EQphi}.
\end{proposition}
\begin{remark} \label{rem:2.1}
We  show below in section~\ref{sec3.2} that the logarithmic potential of the measure $\lambda$ defined in \eqref{EqMesMMI} is equal to the exponent of the main term of the asymptotics of the multiple Meixner polynomials for all values of  parameters   $(c_1,c_2)\in \{(0,1)\times(0,1)\}$ (see Theorem~\ref{T3.2}). 
Correspondently  (for all values of  the parameters) the zero counting measure of the polynomials converges to $\lambda$
$$
\nu_{\widetilde{M}_{2n}^{\beta;\vec{c}}(x)}\overset{*}{\to}\lambda,
$$
and the interval between the origin and the smallest real branch point will be the saturation zone. However, for $(c_1,c_2)\in \{(0,1)\times(0,1)\} \setminus \mathcal{N}$ the measure $\lambda$ is not any more a solution of the equilibrium problem 
\eqref{VmesMMI}, \eqref{EqConMMI}, \eqref{SpropMI} with Nikishin's matrix of interaction.  Without going into details, we just mention that measure $\lambda$ from \eqref{EqMesMMI} for $(c_1,c_2)\in \mathcal{A}$ can be presented as a solution of an  equilibrium problem with Angeleso matrix of interaction (see \cite{17}). For  $(c_1,c_2)\in \mathcal{GN}\setminus \mathcal{N}$ it will be more sophisticated equilibrium problem for the condenser with three plates (see examples in \cite{AptKuUMN}, \cite{AptLys}).
It is interesting to mention, that if we take multiple Laguerre polynomials considered in \cite{Lys}, \cite{LysWil} as the direct continuous analogue of the multiple Miexner polynomials we will see a complete correspondence of these three asymptotical regimes ($\mathcal{N}, \mathcal{A}$ and $\mathcal{GN}\setminus\mathcal{N}$) between the continuous and discrete case.
\end{remark}

\subsection{Potential problem and its solution for multiple Meixner polynomials (second kind)}

Let $r=2$ and $0<\beta_{2}<\beta_{1}<\beta_{2}+1$ in the representation
(\ref{weight-2-kind}) for weights of the discrete orthogonality measures. Now
writing the measures  of orthogonality for multiple Meixner polynomials (second
kind) in form of a Nikishin system, we consider the ratio of the weights
$$
\widehat{u}(x)=\frac{w^{\beta_2}(x)}{w^{\beta_1}(x)}=\sum_{k=0}^{\infty}\frac{\nu_{k}}{x-x_{k}},\quad
\nu_{k}=\frac{1}{\Gamma(\beta_{1}-\beta_{2})}\frac{(\beta_{2}-\beta_{1}+1)_{k}}{k!},\quad
x_{k}=-\beta_{2}-k.
$$
This meromorphic function can be understood as the Cauchy transform of a
discrete measure $u$ with support belonging to $\mathbb{R}_-$. Thus
orthogonality measures of multiple Meixner polynomials (second kind)  form a
Nikishin system generated by two discrete measures: $\mu_1$
(see \eqref{mes-2-kind}) with support in $\mathbb{R}_+$ and $u$ with support in
$\mathbb{R}_-$. Since the mass points for both measures are uniformly
distributed over a subset of $\mathbb{R}_+$ and $\mathbb{R}_-$, respectively, then the limiting zero counting measures in the equilibrium problem will be constrained by $dl(x)=dx$ the Lebesgue measure.

Doing a scaling
$$\widetilde{M}_{2n}^{\beta_1,\beta_2;c}(x):=M_{n,n}^{\beta_1,\beta_2;c}(nx)\,,\qquad
\widetilde{R}_{n}^{(\beta_1,\beta_2;c),j}(z)=R_n^{(\beta_1,\beta_2;c),j}(nz)\,
\quad j=1,2\,,$$ we get the varying weight of orthogonality which defines an
external field on $\mathbb{R}_+$.

Thus, we expect that the zero counting measures associated with the scaled multiple orthogonal polynomial and the functions of the second kind weakly converge
$$
\nu_{\widetilde{M}_{2n}^{\beta_1,\beta_2;c}}\overset{*}{\to}\lambda, \qquad
\qquad
\nu_{\widetilde{R}_{n}^{(\beta_1,\beta_2;c),j}}\overset{*}{\to}\mu\,,\quad
j=1,2,
$$
to measures $(\lambda,\mu):=\overrightarrow{\lambda}$
\begin{equation}\label{VmesMMII}
 \begin{cases}
\,\,\abs{\lambda}=2,\quad \lambda\leq l,\quad\supp(\lambda)\subset\mathbb{R}_{+},\\
\,\, \abs{\mu}=1,\quad\mu\leq l,\quad\supp(\mu)\subset\mathbb{R}_{-}.
\end{cases}
\end{equation}
Now  we write the equilibrium problem for these measures
$\vec{\lambda}=(\lambda,\mu)$ as \eqref{vecEQcons}, where:
\begin{itemize}
 \item Matrix of interaction $D$ in the definition \eqref{eq:1.13} of the
 vector potential $U_j^{\vec{\lambda}}$ has the Nikishin's form \eqref{NikMatr}.

\item Scaling of the weights $w^{\beta_1}$  brings external fields $V_1:=(-\ln
 c)\R{x}$ on $\mathbb{R}_+$.
\end{itemize}
Thus, we have the equilibrium relations:
\begin{equation}\label{EqConMMII}
 \begin{cases}
\,\,W_1\,:=\,2\mathcal{P}^{\lambda}-\mathcal{P}^{\mu} -\ln
 c\,\R{x}&\left\{
\begin{array}{l}\leq\kappa_{1},\quad\mbox{on}\quad\supp(\lambda)\subset\mathbb{R}_{+},\\
\geq\kappa_{1},\quad\mbox{on}\quad\mathbb{R}_{+}\cap\supp(l\setminus\lambda).
\end{array}\right.
\\
\,\,W_2\,:=\,2\mathcal{P}^{\mu}-\mathcal{P}^{\lambda} &\left\{
\begin{array}{l}\leq\kappa_{2},\quad\mbox{on}\quad\supp(\mu)\subset\mathbb{R}_{-},\\
\geq\kappa_{2},\quad\mbox{on}
\quad\mathbb{R}_{-}\cap\supp(l\setminus\mu).
\end{array}\right.
\end{cases}
\end{equation}

Now we look for a solution of the equilibrium problem \eqref{VmesMMII},
\eqref{EqConMMII}. In accordance with our approach (by differentiating) we pass
from \eqref{EqConMMII} to the following relations for the   Cauchy transforms of
 the measures $\lambda$ and $\mu$,
\begin{align}\label{conmMII}
\quad \hat{\lambda}&=\hat{\mu}-\hat{\lambda}-\ln c,&\mbox{on}
\quad \supp(\lambda)\cap\supp(l\setminus\lambda)=:\Sigma_\lambda,\\
\quad \hat{\mu}&=\hat{\lambda}-\hat{\mu},&\mbox{on}\quad
\supp(\mu)\cap\supp(l\setminus\mu)=:\Sigma_\mu.\notag
\end{align}
We point out that if we define on the three sheeted Riemann surface
$$
\mathcal{R} \,:=\, \left(\mathcal{R}_0, \mathcal{R}_1, \mathcal{R}_2
\right)\,:\qquad
\begin{cases}
\,\,\mathcal{R}_0\,:=\,\overline{\mathbb{C}}\setminus \Sigma_\lambda\,,\\
\,\,\mathcal{R}_1\,:=\,\overline{\mathbb{C}}\setminus \left\{\,\Sigma_\lambda\bigcup
\Sigma_\mu\right\}\,,\\
\,\,\mathcal{R}_2\,:=\,\overline{\mathbb{C}}\setminus  \Sigma_\mu\,.
\end{cases}
$$
a function
$$
H(z):=\left\{\begin{array}{ll}
\hat{\lambda},\quad & \mbox{in}\quad \mathcal{R}_0,\\
\hat{\mu}-\hat{\lambda}-\ln c,\quad & \mbox{in}\quad \mathcal{R}_{1},\\
-\hat{\mu}-\ln c,\quad & \mbox{in}\quad \mathcal{R}_{2},\end{array}\right.
$$
then due to \eqref{conmMII}  the branches of this function have analytic
continuation one to another throw the cuts joining the sheets
$\left(\mathcal{R}_0, \mathcal{R}_1, \mathcal{R}_2 \right)$ of the Riemann
surface $\mathcal{R}$. Observe that from \eqref{VmesMMII} follows that $H(z)$ behaves near infinity as
$$
H(z)=\left\{\begin{array}{ll}
\dfrac{2}{z}+\cdots,\quad&\mbox{as $z\to\infty\,\,$ in}\quad \mathcal{R}_{0},\\
-\ln c -\dfrac{1}{z}+\cdots,\quad&\mbox{as $z\to\infty\,\,$ in} \quad
\mathcal{R}_{1}\, \mbox{and}\,\, \mathcal{R}_{2}.\end{array}\right.
$$
In the same way like we did it in the previous subsections  we model the
constrain condition $(\lambda\leq l$ on $\mathbb{R_+}$ and $\mu\leq l$ on
$\mathbb{R_-}$, $dl(x)=dx)$  representing the function $H(z)$ as logarithm of a
meromorphic function $\phi(z)$ on $\mathcal{R}$, i.e.
\begin{equation}\label{HnMII}
H(z)=\ln\phi(z)\,,\qquad \phi \in \mathcal{M}(\mathcal{R}).
\end{equation}
To have a jump for the function $H$ along $\mathbb{R_+}$ lifted to sheet
$\mathcal{R}_{0}$ and along $\mathbb{R_-}$ lifted to sheet $\mathcal{R}_{1}$
we force the function $\phi(z)$ to have pole $\phi=\infty$ of order 2 at the point
$z=0$ on the sheet $\mathcal{R}_{1}$, and zeros $\phi=0$ of order 1 at the
points $z=0$ on the two sheets  $\mathcal{R}_{0}$ and  $\mathcal{R}_{2}$
 \begin{equation}
\phi(z)=\left\{\begin{array}{ll}\,\, 0\,,\quad&\mbox{as $z\to 0\,\,$ in}\quad
\mathcal{R}_{0},\\
\infty^2,\quad&\mbox{as $z\to 0\,\,$ in}\quad \mathcal{R}_{1},\\
\,\,0\,,\quad&\mbox{as $z\to 0\,\,$ in}\quad \mathcal{R}_{2}.\end{array}\right.
\label{behavior 0_phi_meixner2}
\end{equation}
Hence,
$H\in\mathcal{M}\left(\mathcal{R}\setminus\left\{\mathbb{R}_+^{(0)}\bigcup\mathbb{R}^{(1)}\bigcup\mathbb{R}_-^{(2)}\right\}\right)$.
From the behavior of $H$ at infinity we conclude
\begin{equation}
\phi(z)=\left\{\begin{array}{ll} 1+\dfrac{2}{z}+\cdots,\quad&\mbox{as
$z\to\infty\,\,$ in}\quad
\mathcal{R}_{0},\\
\dfrac{1}{ c}\left(1-\dfrac{1}{z}\right)+\cdots,\quad&\mbox{as $z\to\infty\,\,$ in}
\quad \mathcal{R}_{1}\, \mbox{and}\,\, \mathcal{R}_{2}.\end{array}\right.
\label{behavior_infty phi_meixner2}
\end{equation}
Now we obtain an equation for the algebraic
function $\phi$. We write
expansions of  branches $\phi$ up to $\mathcal{O}(1/z^3)$
\begin{eqnarray*}
\phi_0(z) & = & 1+ \frac{2}{z} + \frac{A}{z^2} + \mathcal{O}\left(z^{-3}\right),
\\ \phi_1(z) & = & \sigma - \frac{\sigma}{z} + \frac{B}{z^2} + \mathcal{O}\left(z^{-3}\right),\quad\mbox{where}\,\,\sigma=c^{-1},
 \\
 \phi_2(z) & = & \sigma - \frac{\sigma}{z} + \frac{D}{z^2} + \mathcal{O}\left(z^{-3}\right).
 \end{eqnarray*}
 Since $\phi$ is a rational function on the Riemann surface $\mathcal{R}$
and all its zeros and poles  are shown in \eqref{behavior 0_phi_meixner2}, we
consider the Vieta relations for the product and for the pairwise products of
these expansions and obtain the following linear system for unknown $A$ and
$B+D$:
\begin{equation*}
\left\{\begin{array}{ll} (B+D) + 3\,\sigma+\sigma\,A&=\,0\\
(1+\sigma)\,(B+D) -4\sigma+\sigma^2+2\sigma\,A)&=\,0\end{array}\right.\quad \Rightarrow \quad
A:=\frac{4\sigma-1}{\sigma-1}\,,\quad B+D:=\frac{\sigma^2+2\sigma}{\sigma-1}\,.
\end{equation*}
Hence, we obtain
\begin{equation}\label{phieqII}
\phi^3-\frac{(2x^2-2x-1)\sigma+(x+1)^2}{x^2}\phi^2+\sigma\frac{(x-2)\sigma+2(x+1)}{x}\phi-\sigma^2=0.
\end {equation}
This algebraic function has four branch points, namely infinity point $\infty$ and three real points $e_-, e_1, e_2$ ($e_-<0<e_1<e_2$).  The cut joining the logarithmic singularities of the function $H$ starts at  point $0$ on sheet $\mathcal{R}_{0}$ goes along $\mathbb{R}_+$ to the branch point $e_1$, passing  to  sheet $\mathcal{R}_{1}$ goes there  back via point $0$ to the branch point $e_-$ and then again  passing  to  sheet $\mathcal{R}_{2}$ goes there back to point $0$. Thus we have, indeed
 $$H \in
\mathcal{M}\left(\,\mathcal{R}\,\setminus\,\{\,[0,e_1]^{(0)}\bigcup\,[e_1,e_-]^{(1)}
\bigcup\,[e_-,0]^{(2)}\}\right)\,.$$ Recall that the density of the measure
$\lambda$ is determined by the limiting values at the points $[0,e_2]$ of the
imaginary part of the $\ln\phi(z)/\pi$
 or equivalently,
\begin{equation}\label{EqMesMMII}
\frac{d\lambda(x)}{dx}\,=\,\lim_{y\to 0}\,\frac{1}{\pi}\, \Ima \,H(x+i
y)\,,\qquad x\in [0,\, e_2]\quad \mbox{in}\quad \mathcal{R}_{0}\,.
\end{equation}
Since function $\phi$ is negative on $[0,\, e_1]\quad \mbox{in}\quad
\mathcal{R}_{0}$ we have that density of the measure $\lambda$ is equal unity
there. Therefore, $[0,\, e_1]$ corresponds to the saturation zone and $\Sigma = [e_1,\, e_2]$ is the equilibrium zone. Thus, we have proved the following result.
\begin{proposition} \label{prop:2.2}
The density of the absolute continues measure $\lambda$,  solving equilibrium
problem \eqref{VmesMMII}, \eqref{EqConMMII} is equal to \eqref{EqMesMMII}, where
the function $H$ is defined by \eqref{HnMII} and \eqref{phieqII}.
\end{proposition}

\section{$n$th-root asymptotic from recurrence relations and equilibrium problem }\label{nroot}

In this section, starting from the coefficient of the recurrence relations, we
obtain the main term of asymptotics of multiple Meixner polynomials and then we
check the connection of this term with the equilibrium problem studied in the previous section. The main steps of our approach here are the following.

1) We begin by forming a transition matrix for our recursion. We consider vectors
$\overrightarrow{V}_1$, $\overrightarrow{V}_2$ and $\overrightarrow{V}_3$ with coordinates
\begin{equation}\label{3.1}
\frac{M_{n_1,n_2}(x)}{(n_1+n_2)!}\;,
\end{equation}
taken for the following values of the multi-index $(n_1,n_2)$:
$$
[(n-1,n-1);\;(n,n-1);\;(n,n)]\rightarrow\overrightarrow{V}_1\;,
$$
$$
[(n,n-1);\;(n,n);\;(n+1,n)]\rightarrow\overrightarrow{V}_2\;,
$$
$$
[(n,n);\;(n+1,n);\;(n+1,n+1)]\rightarrow\overrightarrow{V}_3\;.
$$
We define transition matrices $A_n^{(j)},\,j=1,2$, by means of relations:
$$
A_n^{(1)}\,\overrightarrow{V}_1=\overrightarrow{V}_2,\quad
A_n^{(2)}\,\overrightarrow{V}_2=\overrightarrow{V}_3\,.
$$
To obtain the coefficients of these transition matrices we need two types of recurrence relations for $M_{n_1,n_2}$, which connects the following indices
\begin{equation}\label{3.2}
\left\{\begin{array}{l}
I) \quad (n,n)\qquad\mbox{and}\qquad(n+1,n),\;(n,n),\;(n,n-1),\;(n-1,n-1),\\
\\
II) \quad (n+1,n)\qquad\mbox{and}\qquad(n+1,n+1),\;(n+1,n),\;(n,n),\;(n,n-1).
\end{array}\right.\end{equation}
The first type of the recursions we already have. They follow from  (\ref{MIrec}) and (\ref{MIIrec})
if we set $(n_1,n_2)=(n,n)$. To get the recursions of the second type we substitute
$(n_1,n_2)=(n,n+1)$ in  (\ref{MIrec}) and (\ref{MIIrec}) and then swap
$a_1$ with $a_2$ and $\beta_1$ with $\beta_2$, respectively.

2) Analysis of  recursions for polynomials (\ref{3.1}) with indices (\ref{3.2}) (these recursions follows from (\ref{MIrec}) and (\ref{MIIrec})) shows that if we
change the variable there as
$$
t:=\frac{x}{n}\;,
$$
then we get the recurrence relations with coefficients having limits when $n\to\infty$.
It gives us a regime
\begin{equation}\label{limreg}
\left\{\begin{array}{l}
n\to\infty,\\
\\
t\in K\Subset\mathbb{C},
\end{array}\right.\end{equation}
for investigation of asymptotics starting from the recurrence relations.
Thus we form a limiting transition matrix
\begin{equation}\label{2.6}
A(t)=\lim\limits_{n\to\infty}A_n^{(2)}\,A_n^{(1)},\quad t\in K\Subset\mathbb{C},
\end{equation}
and find eigenvalues of this matrix. We fix the following order for these eigenvalues
$$
|L_1(t)|\geq|L_2(t)|\geq|L_3(t)|,\quad t \in \mathbb{C}.
$$

3) Then, $n$th root asymptotics for polynomials \eqref{3.1} can be obtained by means of the following lemma.
\begin{lemma} There holds the following asymptotic formula (when $n \to \infty$)
\begin{equation}\label{3.3}
\frac{1}{n}\ln\left|\frac{M_{n,n}(x)}{(2n)!}\right|=
\frac{1}{n}\int\limits_0^n\ln\left|L_1\left(\frac{x}{\tilde{n}}\right)\right|\,d\tilde{n}+o(1),
\end{equation}
uniformly for $\frac{x}{n}\in K\Subset\Omega$, where $\Omega$ is a domain containing $\infty$-point and bounded by a curve $\Gamma:=\,\{t:L_1(t)=L_2(t)\}$.
\end{lemma}
The proof of this Lemma follows from Poincare's theorem on ratio asymptotics of
solutions of the recurrence relations (see \cite{Poun}) and from the theorem of  Kuijlaars and Van~Assche on the $n$th root asymptotics of the solutions of the recurrence relations with varying coefficients (see \cite{vanAsscheKU} as well as \cite{APTvanAssche} and \cite{AptGERvanAssche}).

We note that function  $L(t)$ is an algebraic function  since it is a root of the polynomials in $L$ with rational coefficients in $t$ (which is the characteristic polynomial of the matrix $A(t)$). This observation allows us to use uniformization of the algebraic curve $L(t)$
$$
\left\{\begin{array}{l} L=F_1(s)\\
n=\displaystyle\frac{x}{F_2(s)}
\end{array}\right.,\qquad \mbox{where $s$ is the parameter of uniformization},
$$
to evaluate explicitly  the integral in the right hand side of \eqref{3.3} .

4) Finally we check that the obtained explicit expression for the right hand side of \eqref{3.3} is the logarithmic potential of measure  $\lambda$,
which is a solution of the corresponding vector equilibrium problem from the previous section
$$
\frac{1}{n}\int\limits_0^n\,\ln\left|L_1\left(\frac{x}{\tilde{n}}\right)\right|\,d\tilde{n}=\mathcal{P}^{\lambda}(t)\;.
$$
For this purpose we differentiate
$$
\frac{d}{dt}\,\left.\left[\frac{1}{n}\int\limits_0^n\,\ln
\left(L_1\left(\frac{x}{\tilde{n}}\right)\right)\,d\tilde{n}\right]\right|_{n=\frac{x}{t}}=:H(t)\;,
$$
and we verify that
$$
H(t)=\ln\phi_0(t)=\int\frac{d\lambda(\xi)}{t-\xi}\;,
$$
where $\phi(t)$ is the algebraic function from the previous section, which
gives an explicit solution of the equilibrium problem.

\subsection{Asymptotics for multiple Meixner polynomials of the second kind}

We start with the polynomials of the second kind since the calculations for obtaining the asymptotics look simplest. Indeed, the following result is valid.

\begin{theorem}\label{T3} The main term of asymptotics for multiple Meixner orthogonal polynomials of the second kind has the form
$$
\frac{1}{n}\ln\left|\displaystyle\frac{M_{n,n}^{\beta_1,\beta_2,c}(x)}{(2n)!}\right|=\mathcal{P}^\lambda\left(\frac{n}{x}\right)+o(1)\;,
$$
and the convergence is uniformly for $\frac{n}{x}\in K\Subset\mathbb{C}$, where at the right hand side stands the logarithmic potential of the equilibrium measure $\lambda$ for the equilibrium problem \eqref{VmesMMII}, \eqref{EqConMMII}. Furthermore, the Cauchy transform
of measure $\lambda$ has the form
$$
\hat{\lambda}(t)=H(t)=\ln\phi_0(t)\;,\quad t\in\mathbb{C}\backslash[0,e_2]\,,
$$
where $\phi_0(t)=1+\frac{2}{t},\;t\to\infty$ is  branch \eqref{behavior_infty
phi_meixner2} of the algebraic function  $\phi$, defined by equation
\eqref{phieqII}, and $e_2$ is the maximal positive branch points of the
function $\phi$.
\end{theorem}

\noindent{\it Proof.} Following the step 1) of our approach (see description at the beginning of this section) we compute the transition matrices  $A_n^{(1)}$ and $A_n^{(2)}$, and we consider them  in the limit regime  \eqref{limreg}. Doing it we note that the obtained limits coincide. This circumstance allows us to use the matrix
$$
A^{(1)}(t)=\lim\limits_{n\to\infty}A_n^{(1)}\;,\quad t\in K\Subset\mathbb{C}\;,
$$
instead of the transition matrix $A(t)$ from \eqref{2.6}. We just have to multiply by 2 the right hand side of formula \eqref{3.3}.
Hence, we have
$$
A^{(1)}(t)=\left(\begin{array}{ccc}
\frac{1}{2}t-\frac{3}{2}a-1 & -\frac{3}{4}a(a+1) & -\frac{1}{8}a^2(a+1)\\
\\
1 & 0 & 0\\
\\
0 & 1 & 0\\
\end{array}\right)\;.
$$
The characteristic polynomial $P(L,t)$ of the matrix $A^{(1)}(t)$ has the form:
$$
P=L^3-\frac{1}{2}\,L^2t+\frac{3}{2}\,L^2a+\frac{3}{4}a^2L+\frac{3}{4}aL+\frac{1}{8}\,a^3+\frac{1}{8}\,a^2\;.
$$
The equation $P(L,t)=0$ defines an algebraic function  $L(t)$ of genus zero,
therefore this function allows a rational uniformization. Moreover, taking into
account the linear appearance of $t$ in $P(L,t)$ we immediately obtain  this
uniformization if we take
 $L$ as the uniformization parameter:
\begin{equation*}\label{3.9}
\left\{\begin{array}{l} L=L,\\
\\
t=\displaystyle\frac{P(L,t)+\frac{1}{2}L^2t}{\frac{1}{2}L^2}=\frac{1}{4}\,\frac{(2L+b+b^2)\,(2L-b+b^2)\,(2L-1+b^2)}{L^2},
\end{array}\right.
\end{equation*}
where $a$ is changed ($a:=b^2-1$). From here, since  $t=\frac{n}{x}$, we have
\begin{equation}\label{3.10}
n=\frac{4L^2x}{(2L+b+b^2)\,(2L-b+b^2)\,(2L-1+b^2)}\;.
\end{equation}
Now we can integrate by parts the right hand side of \eqref{3.3}, i.e.
$$
2\left(\int\ln(L)\,dn\right)=2n\ln(L)+x\ln(2L+b+b^2)+x\ln(2L-b+b^2)-2x\ln(2L-1+b^2)\;.
$$
Then, we substitute $x(n)$ from \eqref{3.10} into the obtained expression and multiplying it by $1/n$, as a result we obtain the right hand side of \eqref{3.3}
$$
\frac{2}{n}\int\ln(L)\,dn=2\ln(L)+
$$
$$
+\frac{1}{4}\,\frac{1}{L^2}\left((2L+b+b^2)\,(2L-1+b^2)\,\ln\left(\displaystyle\frac{(2L+b+b^2)\,
(2L-b+b^2)}{(2L-1+b^2)^2}\right)\,(2L-b+b^2)\right)=:F(L)\;.
$$
It remains to check, that the real part of $F(L)$ is indeed the logarithmic
potential of the equilibrium measure $\lambda$ of the problem
\eqref{VmesMMII}, \eqref{EqConMMII}. It will be true if the derivative of $F(L)$
with respect to $t$ coincides with the function $H(t)=\ln\phi(t)$, where $\phi(t)$ is the algebraic function \eqref{phieqII}. To differentiate $F$ we use the uniformization variable $L$. Thus,
$$
\frac{d}{dt}F(L(t))=\displaystyle\frac{\frac{d}{dL}F(L)}{\frac{d}{dL}t(L)}=\ln\left(\displaystyle\frac{(2L+b+b^2)\,
(2L-b+b^2)}{(2L-1+b^2)^2}\right)=:H\;.
$$
Finally, substituting
$$
\begin{array}{l}
\phi:=e^H=\displaystyle\frac{(2L+b+b^2)\, (2L-b+b^2)}{(2L-1+b^2)^2}\;;\\
\\
z:=t=\displaystyle\frac{1}{4}\displaystyle\frac{(2L+b+b^2)\,(2L-b+b^2)\,(2L-1+b^2)^2}{L^2}\;;\\
\\
b:=\displaystyle\frac{1}{\sqrt{1-1/\sigma}}\;,\quad \sigma=c^{-1},
\end{array}
$$
we get an equation for the function  $\phi(z)$, i.e.
$$
z^2\phi^3-\phi^2\left(2\sigma z^2+z^2-2z\sigma+2z-\sigma+1\right)+z\sigma(z\sigma+2z-2\sigma+2)\phi-\sigma^2z^2=0\;,
$$
which coincide with \eqref{phieqII}. Theorem is completely proved.\qed

\subsection{Asymptotics for multiple Meixner polynomials of the first kind}\label{sec3.2}

Finally, we consider the Meixner multiple orthogonal polynomials of the first kind.

\begin{theorem}\label{T3.2} The main term of asymptotics for the multiple Meixner orthogonal polynomials of first kind has the form
$$
\frac{1}{n}\ln\left|\displaystyle\frac{M_{n,n}^{\beta;\,\vec{c}}(x)}{(2n)!}\right|=
\mathcal{P}^\lambda\left(\frac{n}{x}\right)\,+\,o(1)\;,
$$
being the convergence uniformly for $\frac{n}{x}\in
K\Subset\mathbb{C}$, where at the right hand side stands the
logarithmic potential of a  measure $\lambda$ defined by its Cauchy
transform which has the form
$$
\hat{\lambda}(t)=H(t)=\ln\phi_0(t)\;,\quad t\in\mathbb{C}\backslash[0,e_2]\,,
$$
where $\phi_0(t)=1+\frac{2}{t},\;t\to\infty$ is  branch \eqref{behavior_infty
phi_meixner2} of the algebraic function  $\phi$, defined by equation
\eqref{EQphi}, and $e_2$ is the maximal positive branch points of the function
$\phi$.

For parameters $(c_1,c_2)\in\mathcal{N}$ (see \eqref{NikDom}) the
measure $\lambda$ is the solution of the equilibrium problem
\eqref{VmesMMI}, \eqref{EqConMMI}, \eqref{SpropMI}.
\end{theorem}

\noindent{\it Proof.} Firstly, we compute the transition matrices  $A_n^{(1)}$ and
$A_n^{(2)}$, and we consider them  in the limit regime \eqref{limreg}
$$
A^{(j)}(t)=\lim\limits_{n\to\infty}A_n^{(j)}\;,\quad t:=\frac{x}{n}\in
K\Subset\mathbb{C}\;.
$$
For the case of the Meixner multiple orthogonal polynomials of first kind these
limiting matrices are different. We have
$$
A^{(1)}=\left(\begin{array}{ccc}
\frac{1}{2}t-\frac{3}{2}a_1-\frac{1}{2}a_2-1 & -\frac{1}{2}a_1(a_1+1)-\frac{1}{2}a_2(a_2+1) & -\frac{1}{2}a_1(a_1+1)(a_1-a_2)\\
\\
1 & 0 & 0\\
\\
0 & 1 & 0\\
\end{array}\right),
$$
and
$$
A^{(2)}=\left(\begin{array}{ccc}
\frac{1}{2}t-\frac{3}{2}a_2-\frac{1}{2}a_1-1 & -\frac{1}{2}a_1(a_1+1)-\frac{1}{2}a_2(a_2+1) & -\frac{1}{2}a_2(a_2+1)(a_2-a_1)\\
\\
1 & 0 & 0\\
\\
0 & 1 & 0\\
\end{array}\right)\;.
$$
Then we compute the characteristic polynomial $P(L,t)$ of the transition matrix
-product $A^{(2)}A^{(1)}(t)$. Like in theorem \ref{T3} the equation $P(L,t)=0$ defines an
algebraic function  $L(t)$ of genus zero. Rational uniformization of this
algebraic curve has the following form

\begin{equation*}
\left\{\begin{array}{l}
L=\displaystyle\frac{(s-a_2(a_2+1)(a_1-a_2))\,(s+a_1(a_1+1)(a_1-a_2))}{(a_1(a_1+1)+(a_2(a_2+1))^2},
\\
\\t=\displaystyle\frac{((s+a_1^2(a_1+1)+a_2^2(a_2+1))\,(s+a_1(a_1+1)^2+a_2(a_2+1)^2)\,(2s+(a_2+a_1+1)(a_1-a_2)^2))}
{((a_1(a_1+1)+a_2(a_2+1))\,(s+a_1(a_1+1)(a_1-a_2))\,(s-a_2(a_2+1)(a_1-a_2)))},
\end{array}\right.
\end{equation*}
with $s$ as a parameter of the uniformization.  From here, since $n =t x$, we
integrate by parts the right hand side of \eqref{3.3}
$$
\int\ln(L)\,dn=n\ln(L)-\ln\left(\displaystyle\frac{s+a_1^2+a_1^3+a_2^2+a_2^3}
{s+a_1+2a_1^2+a_1^3+a_2+2a_2^2+a_2^3}\right)x\,.
$$
Hence, multiplying it by $\frac{1}{n}$ we obtain the right hand side of
\eqref{3.3}

$$
\frac{1}{n}\int\ln(L)\,dn=\ln(L(s))\,+\,
\ln\left(\displaystyle\frac{s+a_1^2+a_1^3+a_2^2+a_2^3}
{s+a_1+2a_1^2+a_1^3+a_2+2a_2^2+a_2^3}\right)t(s)\,.=:F(s)\;.
$$
It remains to check, that the real part of $F$ is indeed the logarithmic
potential of the equilibrium measure  $\lambda$ of the problem \eqref{VmesMMI},
\eqref{EqConMMI}, \eqref{SpropMI}. Thus, we differentiate $F$ using the
uniformization variable $s$. Accordingly, we have
$$
H:=\ln(\phi):=\frac{d}{dt}F=\displaystyle\frac{\frac{d}{ds}F}{\frac{d}{ds}t}=
-\ln\left(\displaystyle\frac{s+a_1^2+a_1^3+a_2^2+a_2^3}{s+a_1+2a_1^2+a_1^3+a_2+2a_2^2+a_2^3}\right)
=:\ln(r(s)).
$$
From here we get a rational uniformization of the curve
$$\phi(z)\,: \qquad \phi\,=\,r(s)\,,\quad z\,=\,t(s),$$
which bring us an equation for the function  $\phi(z)$:
\begin{eqnarray*}
a_2a_1z\phi^3+(-3a_2a_1z-za_1+a_1-za_2+a_2)\phi^2
+(-2+z-a_2-a_1+2za_2+2za_1+3a_2a_1z)\phi\\
-z(a_2+1)(a_1+1)=0,
\end{eqnarray*}
or equivalently
$$
z-\frac{2}{\phi-1}+\frac{a_1+1}{a_1\phi-1-a_1}+\frac{a_2+1}{a_2\phi-a_2-1}=0,
$$
which after substitution $a_1=\displaystyle\frac{c_1}{1-c_1}$ and $a_2=\displaystyle\frac{c_2}{1-c_2}$ coincides with \eqref{EQphi}. 
It remains to apply Proposition~\ref{prop:2.1}. The theorem is proved.\qed

\section{Conclusions and future directions}

In this paper we have investigated the $n$th root asymptotic behavior of the multiple Meixner polynomials of the first and second kind, respectively (see theorems \ref{T3}-\ref{T3.2}).
We have developed an approach based on an algebraic function formulation in connection with some available techniques from logarithmic potential theory. Indeed, for each polynomial sequence we have posed the corresponding equilibrium problem for the limiting zero counting measure and then we solved it in a quite direct way. A detailed analysis of the asymptotic distribution of zeros of multiple Meixner polynomials was given, i.e. we looked at different regions in the complex plane for identifying the support of the extremal measure, the saturation zone, the equilibrium zone, and the free zone.

Finally, other systems of discrete multiple orthogonal polynomials (for $r=2$) might be asymptotically analyzed using the approach presented here
(up to some adaptations) to produce similar outcomes.\\

\noindent\textbf{Acknowledgments.} The authors are grateful to their colleagues Van Assche, Sorokin and Tulyakov for useful discussions.

\parbox{3in}{
\texttt{Jorge Arves\'u Carballo  \\
Department of Mathematics \\
Universidad Carlos III de Madrid \\
Avda. de la Universidad 30 \\
E-28911 Legan\'es, Madrid \\
SPAIN}} \hfil
\parbox{2.5in}{
\texttt{Alexander Aptekarev \\
M. V. Keldysh Institute for Applied Mathematics\\
Russian Academy of Sciences\\
125047 Moscow, Miusskaya pl. 4\\
RUSSIA}}

\end{document}